\documentclass[graybox]{svmult}

\usepackage{mathptmx}       
\usepackage{helvet}         
\usepackage{courier}        
\usepackage{type1cm}        
\usepackage{makeidx}         
\usepackage{graphicx}        
\usepackage{multicol}        
\usepackage[bottom]{footmisc}

\usepackage{mystyle}

\newcommand{\myparagraph}[1]{\vspace{2mm}\noindent\textbf{#1}~}

\makeindex             

\begin{document}

\title*{Low Complexity Regularization \\ of Linear Inverse Problems}
\author{Samuel Vaiter, Gabriel Peyr\'e and Jalal Fadili}
\institute{Samuel Vaiter \at Ceremade, Universit\'e Paris-Dauphine, \email{samuel.vaiter@ceremade.dauphine.fr}
\and 
Gabriel Peyr\'e \at CNRS and Ceremade, Universit\'e Paris-Dauphine, \email{gabriel.peyre@ceremade.dauphine.fr}
\and
Jalal Fadili \at GREYC, CNRS-ENSICAEN-Universit\'e de Caen, \email{Jalal.Fadili@greyc.ensicaen.fr} 
}
%
%
\maketitle




\abstract{Inverse problems and regularization theory is a central theme in imaging sciences, statistics and machine learning. The goal is to reconstruct an unknown vector from partial indirect, and possibly noisy, measurements of it.
A now standard method for recovering the unknown vector is to solve a convex optimization problem that enforces some prior knowledge about its structure. 
This chapter delivers a review of recent advances in the field where the regularization prior promotes solutions conforming to some notion of simplicity/low-complexity. These priors encompass as popular examples sparsity and group sparsity (to capture the compressibility of natural signals and images), total variation and analysis sparsity (to promote piecewise regularity), and low-rank (as natural extension of sparsity to matrix-valued data). 
Our aim is to provide a unified treatment of all these regularizations under a single umbrella, namely the theory of partial smoothness. This framework is very general and accommodates all low-complexity regularizers just mentioned, as well as many others. Partial smoothness turns out to be the canonical way to encode low-dimensional models that can be linear spaces or more general smooth manifolds.
This review is intended to serve as a one stop shop toward the understanding of the theoretical properties of the so-regularized solutions. It covers a large spectrum including: (i) recovery guarantees and stability to noise, both in terms of $\ell^2$-stability and model (manifold) identification; (ii) sensitivity analysis to perturbations of the parameters involved (in particular the observations), with applications to unbiased risk estimation ; (iii) convergence properties of the forward-backward proximal splitting scheme, that is particularly well suited to solve the corresponding large-scale regularized optimization problem.}

\section{Inverse Problems and Regularization}
\label{sec:intro-context}

In this chapter, we deal with finite-dimensional linear inverse problems.

\subsection{Forward Model}
\label{sec:intro-forward}

Let $x_0 \in \RR^\N$ be the unknown vector of interest. Suppose that we observe a vector $y \in \RR^\P$ of $\P$ linear measurements according to
\eql{\label{eq:lin-inverse-problem}
  y = \Phi x_0 + w ,
}
where $w \in \RR^\P$ is a vector of unknown errors contaminating the observations. The forward model \eqref{eq:lin-inverse-problem} offers a model for data acquisition that describes a wide range of problems in data processing, including signal and image processing, statistics, and machine learning. The linear operator $\Phi : \RR^\N \rightarrow \RR^\P$, assumed to be known, is typically an idealization of the acquisition hardware in imaging science applications, or the design matrix in a parametric statistical regression problem. The noise $w$ can be either deterministic (in this case, one typically assumes to know some bound on its $\ell^2$ norm $\norm{w}$), or random (in which case its distribution is assumed to be known). Except in Sections~\ref{sec-model-consistency-proba} and~\ref{sec-risk-estimation} where the noise is explicitly assumed random, $w$ is deterministic throughout the rest of the chapter. We refer to~\cite{refregier2004statistical} and~\cite{bovik2005noise} for a comprehensive account on noise models in imaging systems.

Solving an inverse problem amounts to recovering $x_0$, to a good approximation, knowing $y$ and $\Phi$ according to~\eqref{eq:lin-inverse-problem}. Unfortunately, the number of measurements $\P$ can be much smaller than the ambient dimension $\N$ of the signal. Even when $\P = \N$, the mapping $\Phi$ is in general ill-conditioned or even singular. This entails that the inverse problem is in general ill-posed. In signal or image processing, one might for instance think of $\Phi$ as a convolution with the camera point-spread function, or a subsampling accounting for low-resolution or damaged sensors. In medical imaging, typical operators represent a (possibly subsampled) Radon transform (for computerized tomography), a partial Fourier transform (for magnetic resonance imaging), a propagation of the voltage/magnetic field from teh dipoles to the sensors (for electro- or magneto-encephalography). In seismic imaging, the action of $\Phi$ amounts to a convolution with a wavelet-like impulse response that approximates the solution of a wave propagation equation in media with discontinuities. For regression problems in statistics and machine learning, $\Phi$ is the design matrix whose columns are $\P$ covariate vectors.

\subsection{Variational Regularization}
\label{sec:intro-reg}

As argued above, solving an inverse problem from the observations~\eqref{eq:lin-inverse-problem} is in general ill-posed. In order to reach the land of well-posedness, it is necessary to restrict the inversion process to a well-chosen subset of $\RR^\N$ containing the plausible solutions including $x_0$; e.g. a linear space or a union of subspaces. A closely related procedure, that we describe next, amounts to adopting a variational framework where the sought-after solutions are those where a prior penalty/regularization function is the smallest. Though this approach may have a maximum a posteriori Bayesian interpretation, where a random prior is placed on $x_0$, this is not the only interpretation. In fact, we put no randomness whatsoever on the class of signals we look for. We will not elaborate more on these differences in this chapter, but the reader may refer to~\cite{gribonval2011should} for an insightful discussion.



The foundations of regularization theory can be traced back to the pioneering work of the Russian school, and in particular of Tikhonov in 1943 when he proposed the notion of conditional well-posedness. In 1963, Tikhonov \cite{Tikhonov63_1,Tikhonov63_2} introduced what is now commonly referred to as Tikhonov (or also Tikhonov-Phillips) regularization, see also the book~\cite{tikhonov1977solutions}. This corresponds, for $\la>0$, to solving an optimization problem of the form
\eql{\label{eq-lagrangian}\tag{$\Pp_{y,\lambda}$}
  \xsol \in \uArgmin{x \in \RR^\N} \frac{1}{2\la} \norm{\Phi x - y}^2 +  J(x).
}

\subsubsection{Data fidelity} 

In~\eqref{eq-lagrangian}, $\norm{\Phi x - y}^2$ stands for the data fidelity term. If the noise happens to be random, then using a likelihood argument, an appropriate fidelity term conforming to the noise distribution can be used instead of the quadratic data fidelity. Clearly, it is sufficient then to replace the latter by the negative log-likelihood of the distribution underlying the noise. Think for instance of the Csisz\'ar's I-divergence for Poisson noise.
We would also like to stress that many of the results provided in this chapter extend readily when the quadratic loss in the fidelity term, i.e. $\mu \mapsto \norm{y - \mu}^2$, is replaced by any smooth and strongly convex function, see in particular Remark~\ref{rem-generic-loss}.
To make our exposition concrete and digestible, we focus in the sequel on the quadratic loss. 

\subsubsection{Regularization} 

The function $J : \RR^\N \to \RR$ is the regularization term which is intended to promote some prior on the vector to recover. We will consider throughout this chapter that $\J$ is a convex finite-valued function. Convexity plays an important role at many locations, both on the recovery guarantees and the algorithmic part. 
See for instance Section~\ref{sec-algorithms} which gives a brief overview of recent algorithms that are able to tackle this class of convex optimization problems. 
It is however important to realize that non-convex regularizing penalties, as well as non-variational methods (e.g. greedy algorithms), are routinely used for many problems such as sparse or low-rank recovery. They may even outperform in practice their convex counterparts/relaxation. It is however beyond the scope of this chapter to describe these algorithms and the associated theoretical performance guarantees. We refer to Section~\ref{sec-model-selection} for a brief account on non-convex model selection approaches.

The scalar $\lambda>0$ is the regularization parameter. It balances the trade-off between fidelity and regularization. Intuitively, and anticipating on our theoretical results hereafter, this parameter should be adapted to the noise level $\norm{w}$ and the known properties of the vector $x_0$ to recover. Selecting optimally and automatically $\la$ for a given problem is however difficult in general. This is at the hear of Section~\ref{sec:intro-sensitivity}, where unbiased risk estimation strategies are shown to offer a versatile solution. 

Note that since $\Phi$ is generally not injective and $\J$ is not coercive, the objective function of~\eqref{eq-lagrangian} is neither coercive nor strictly convex. In turn, there might be existence (of minimizers) issues, and even if minimizers exist, there are not unique in general.

Under mild assumptions, problem~\eqref{eq-lagrangian} is formally equivalent to the constrained formulations
\begin{align}
	& \tag{$\Pp_{y,\epsilon}^1$} \min \enscond{ J(x) }{ \norm{y - \Phi x} \leq \epsilon }, \label{eq-constraint-epsilon} \\
	& \tag{$\Pp_{y,\gamma}^2$} \min \enscond{ \norm{y-\Phi x} }{ J(x) \leq \gamma }, \label{eq-constraint-gamma}
\end{align}
in the sense that there exists a bijection between each pair of parameters among $(\la,\epsilon,\gamma)$ so that the corresponding problems share the same set of solutions. However, this bijection is not explicit and depends on $y$, so that both from an algorithmic point of view and a theoretical one, each problem may need to be addressed separately. See the recent paper \cite{Ciak13} and references therein for a detailed discussion, and~\cite[Theorem~2.3]{Lorenz-Necessary} valid also in the non-convex case. We focus in this chapter on the penalized/Tikhonov formulation~\eqref{eq-lagrangian}, though most of the results stated can be extended to deal with the constrained ones~\eqref{eq-constraint-epsilon} and~\eqref{eq-constraint-gamma} (the former is known as the residual method or Mozorov regularization and the latter as Ivanov regularization in the inverse problems literature).

The value of $\la$ should typically be an increasing function of $\norm{w}$. In the special case where there is no noise, i.e. $w=0$, the fidelity to data should be perfect, which corresponds to considering the limit of \eqref{eq-lagrangian} as $\la \rightarrow 0^+$. Thus, assuming that $y \in \Im(\Phi)$, as is the case when $w=0$, it can be proved that the solutions of~\eqref{eq-lagrangian} converge to the solutions of the following constrained problem \cite{Tikhonov63_1,scherzer2009variational}
\eql{\label{eq:reg-noiseless}\tag{$\Pp_{y,0}$}
  \xsol \in \uArgmin{x \in \RR^\N} J(x) 
  \qsubjq
  \Phi x = y .
}

\subsection{Notations}

For any subspace $T$ of $\RR^\N$, we denote $\proj_T$ the orthogonal projection onto $T$, $x_T = \proj_T(x)$ and $\Phi_T = \Phi \proj_T$. For a matrix $A$, we denote $A^*$ its transpose, and $A^+$ its Moore-Penrose pseudo-inverse. For a convex set $E$, $\aff(E)$ denotes its affine hull (i.e. the smallest affine space containing it), and $\Lin(E)$ its linear hull (i.e. the linear space parallel to $\aff(E)$). Its relative interior $\ri(E)$ is the interior for the topology of $\aff(E)$ and $\rbd(E)$ is its relative boundary. For a manifold $\Mm$, we denote $\tgtManif{\x}{\Mm}$ the tangent space of $\Mm$ at $x \in \Mm$. A good source on smooth manifold theory is \cite{lee2003smooth}. 

A function $J: \RR^\N \rightarrow \RR \cup \{+\infty\}$ is said to be proper if it is not identically $+\infty$. It is said to be finite-valued if $\J(x) \in \RR$ for all $x \in \RR^\N$. We denote $\dom(J)$ the set of points $x$ where $J(x) \in \RR$ is finite. $J$ is said to be closed if its epigraph $\enscond{(x,y)}{J(x) \leq y}$ is closed. For a set $C \subset \RR^N$, the indicator function $\iota_C$ is defined as $\iota_C(x)=0$ if $x \in C$ and $\iota_C(x)=+\infty$ otherwise.

We recall that the subdifferential at $x$ of a proper and closed convex function $J : \RR^\N \rightarrow \RR \cup \{+\infty\}$ is the set
\eq{
	\partial J(x) = \enscond{\eta \in \RR^N}{ \foralls \delta \in \RR^\N, \;
		J(x+\delta) \geq J(x) + \dotp{\eta}{\delta} }.
}
Geometrically, when $J$ is finite at $x$, $\partial J(x)$ is the set of normals to the hyper-planes supporting the graph of $J$ and tangent to it at $x$. Thus, $\partial J(x)$ is a closed convex set. It is moreover bounded, hence compact, if and only if $x \in \Int(\dom(J))$. 
The size of the subdifferential at $x \in \dom(J)$ reflects in some sense the degree of non-smoothness of $J$ at $x$. The larger the subdifferential at $x$, the larger the ``kink'' of the graph of $J$ at $x$. 
In particular, if $J$ is differentiable at $x$, then $\partial J(x)$ is a singleton and $\partial J(x) = \{ \nabla J(x) \}$.

As an illustrative example, the subdifferential of the absolute value is
\eql{\label{eq-subdiff-absval}
	\foralls x \in \RR, \quad
	\partial |\cdot|(x) = \choice{
		\sign(x) \qifq x \neq 0, \\
		{[-1,1]} \quad \text{otherwise.}
	}
}
The $\ell^1$ norm
\eq{
	\foralls x \in \RR^N, \quad \norm{x}_1 = \sum_{i=1}^N |x_i|
}
is a popular low-complexity prior (see Section~\ref{subsec-l1example} for more details). Formula~\eqref{eq-subdiff-absval} is extended by separability to obtain the subdifferential of the $\ell^1$ norm
\eql{\label{eq-formula-l1-subdiff}
	\partial\norm{\cdot}_1(x) = \enscond{ \eta \in \RR^N }{ \norm{\eta}_\infty \leq 1 
		\qandq \foralls i \in I, \; \sign(\eta_i)=\sign(x_i)  }
}
where $I=\supp(x)=\enscond{i}{x_i \neq 0}$. Note that at a point $x \in \RR^N$ such that $x_i \neq 0$ for all $i$, $\norm{\cdot}_1$ is differentiable, and $\partial\norm{\cdot}_1(x) = \{ \sign(x) \}$.



\section{Low Complexity Priors}
\label{sec:intro-low}

A recent trend in signal and image processing, statistics and machine learning is to make use of large collections of so-called ``models'' to account for the complicated structures of the data to handle. Generally speaking, these are manifolds $\Mm$ (most of the time linear subspaces), and hopefully of low complexity (to be detailed later), that capture the properties of the sought after signal, image or higher dimensional data. In order to tractably manipulate these collections, the key idea underlying this approach is to encode these manifolds in the non-smooth parts of the regularizer $J$. As we detail here, the theory of partial smoothness turns out to be natural to provide a mathematically grounded and unified description of these regularizing functions.

\subsection{Model Selection}
\label{sec-model-selection}

The general idea is thus to describe the data to recover using a large collection of models $\mathbb{M} = \{\Mm\}_{\Mm \in \mathbb{M}}$, which are manifolds. 
The ``complexity'' of elements in such a manifold $\Mm$ is measured through a penalty $\pen(\Mm)$. A typical example is simply the dimensionality of $\Mm$, and it should reflect the intuitive notion of the number of parameters underlying the description of the vector $x_0 \in \Mm$ that one aims at recovering from the noisy measurements of the form~\eqref{eq:lin-inverse-problem}. As popular examples of such low complexity, one thinks of sparsity, piecewise regularity, or low rank.
%
%
Penalizing in accordance to some notion of complexity is a key idea, whose roots can be traced back to the statistical and information theory literature, see for instance~\cite{mallows,akaike1973information}. 

Within this setting, the inverse problem associated to the measurements~\eqref{eq:lin-inverse-problem} is solved by restricting the inversion to an optimal manifold as selected by $\pen(\Mm)$. Formally, this would correspond to solving~\eqref{eq-lagrangian} with the combinatorial regularizer
\eql{\label{eq:intro-selec-pen}
  	J(x) = \inf \enscond{ \pen(\Mm) }{ \Mm \in \mathbb{M} \qandq x \in \Mm  }.
}

A typical example of such a model selection framework is that with sparse signals, where the collection $\mathbb{M}$ corresponds to a union of subspaces, each of the form 
\eq{
	\Mm = \enscond{x \in \RR^\N}{\supp(x) \subseteq I}.
}
Here $I \subseteq \{1,\ldots,N\}$ indexes the supports of signals in $\Mm$, and can be arbitrary. In this case, one uses $\pen(\Mm) = \dim(\Mm) = |I|$, so that the associated combinatorial penalty is the so-called $\ell^0$ pseudo-norm
\eql{\label{eq-lzero-pseudonorm}
	J(x) = \norm{x}_0 = \abs{ \supp (x) }  = |\enscond{i \in \{1,\ldots,\N\}}{ x_i \neq 0 }|.
}
Thus, solving~\eqref{eq-lagrangian} is intended to select a few active variables (corresponding to non-zero coefficients) in the recovered vector.

These sparse models can be extended in many ways. For instance, piecewise regular signals or images can be modeled using manifolds $\Mm$ that are parameterized by the locations of the singularities and some low-order polynomial between these singularities. The dimension of $\Mm$ thus grows with the number of singularities, hence the complexity of the model. 

\myparagraph{Literature review.}
The model selection literature~\cite{birge1997model,barron1999risk,birge2007minimal} proposes many theoretical results to quantify the performance of these approaches. However, a major bottleneck of this class of methods is that the corresponding $J$ function defined in~\eqref{eq:intro-selec-pen} is non-convex, and even not necessarily closed, thus typically leading to highly intractable combinatorial optimization problems. For instance, in the case of the $\ell^0$ penalty~\eqref{eq-lzero-pseudonorm} and for an arbitrary operator $\Phi$, \eqref{eq-lagrangian} is known to be NP-hard, see e.g.~\cite{natarajan1995sparse}.

It then appears crucial to propose alternative strategies which allow to deploy fast computational algorithms. A first line of work consists in finding stationary points of~\eqref{eq-lagrangian} using descent-like schemes. For instance, in the case of the $\ell^0$ pseudo-norm, this can be achieved using iterative hard thresholding~\cite{BlumensathHardThresh,Attouch13}, or iterative reweighting schemes which consist of solving a sequence of weighted $\ell^1$- or $\ell^2$-minimization problems where the weights used for the next iteration are computed from the values of the current solution, see for instance~\cite{RaoFOCUSS,CandesReweighting,DaubechiesIRLS} and references therein. 
Another class of approaches is that of greedy algorithms. These are algorithms which explore the set of possible manifolds $\Mm$ by progressively, actually in a greedy fashion, increasing the value of $\pen(\Mm)$. The most popular schemes are matching pursuit~\cite{mallat1993matching} and its orthogonal variant~\cite{pati1993orthogonal,davis94adaptive}, see also the comprehensive review~\cite{needell2008greedy} and references therein.
The last line of research, which is the backbone of this chapter, consists in considering convex regularizers which are built in such away that they promote the same set of low-complexity manifolds $\mathbb{M}$. In some cases, the convex regularizer proves to be the convex hull of the initial (restricted) non-convex combinatorial penalty~\eqref{eq:intro-selec-pen}. But these convex penalties can also be designed without being necessarily convexified surrogates of the original non-convex ones. 

In the remainder of this section, we describe in detail a general framework that allows model selection through the general class of convex partly smooth functions.

\subsection{Encoding Models into Partly Smooth Functions}
\label{sec-partial-smoothness}

Before giving the precise definition of our class of convex priors, we define formally the subspace $T_x$.

\begin{definition}[Model tangent subspace]
  For any vector $x \in \RR^N$, we define the \emph{model tangent subspace} of $x$ associated to $\J$
  \begin{equation*}
    \T_x = \Lin (\partial \J(x))^\bot, 
  \end{equation*}
\end{definition}
In fact, the terminology ``tangent'' originates from the sharpness property of Definition~\ref{dfn-partly-smooth}\eqref{PS-Sharp} below, when $x$ belongs to the manifold $\Mm$.

When $J$ is differentiable at $x$, i.e. $\partial J(x)=\ens{\nabla J(x)}$, one has $\T_x=\RR^N$. On the contrary, when $\J$ is not smooth at $x$, the dimension of $\T_x$ is of a strictly smaller dimension, and $\J$ essentially promotes elements living on or close to the affine space $x+T_x$. 

We can illustrate this using the $\ell^1$ norm $J=\norm{\cdot}_1$ defined in~\eqref{eq-subdiff-absval}. Using formula~\eqref{eq-formula-l1-subdiff} for the subdifferential, one obtains that
\eq{
	T_{\x} = \enscond{ u \in \RR^\N }{ \supp(u) \subseteq \supp(\x) }, 
}
which is the set of vector having the same sparsity pattern as $\x$. 

Toward the goal of studying the recovery guarantees of problem~\eqref{eq:intro-selec-pen}, our central assumption is that $J$ is a partly smooth function relative to some manifold $\Mm$. Partial smoothness of functions was originally defined~\cite{lewis2002active}. Loosely speaking, a partly smooth function behaves smoothly as we move on the manifold $\Mm$, and sharply if we move normal to it. Our definition hereafter specializes that of~\cite{lewis2002active}  to the case of finite-valued convex functions.

\begin{definition} \label{dfn-partly-smooth}
	Let $J$ be a finite-valued convex function. 
	$J$ is \emph{partly smooth at $\x$ relative to a set $\Mm$} containing $\x$ if 
	\begin{enumerate}[(i)] \setlength{\itemsep}{0pt}
		\item\label{PS-C2}(Smoothness) $\Mm$ is a $C^2$-manifold around~$\x$ and $J$ restricted to $\Mm$ is $C^2$ around~$\x$.
		\item\label{PS-Sharp}(Sharpness) The tangent space $\tgtManif{\x}{\Mm}$ is $T_{\x}$.
		\item\label{PS-DiffCont}(Continuity) The set-valued mapping $\partial J$ is continuous at $\x$ relative to~$\Mm$.
	\end{enumerate}
	$J$ is said to be \emph{partly smooth relative to a set $\Mm$} if $\Mm$ is a manifold and $J$ is partly smooth at each point $\x \in \Mm$ relative to $\Mm$.
    $J$ is said to be \emph{locally partly smooth at $\x$ relative to a set $\Mm$} if $\Mm$ is a manifold and there exists a neighbourhood $U$ of $\x$ such that $J$ is partly smooth at each point of $\Mm \cap U$ relative to $\Mm$.
\end{definition}

\begin{remark}[Uniqueness of $\Mm$]
In the previous definition, $\Mm$ needs only to be defined locally around $\x$, and it can be shown to be locally unique, see~\cite[Corollary~ 4.2]{hare2007identif}. In the following we will thus often denote $\Mm_x$ any such a manifold for which $J$ is partly smooth at $x$.  
\end{remark}

Taking once again the example of $J=\norm{\cdot}_1$, one sees that in this case, $\Mm_x=T_x$ because this function is polyhedral. Section~\ref{sec-spectral-priors} bellow defines fonctions $J$ for which $\Mm_x$ differs in general from $T_x$.

\subsection{Examples of Partly Smooth Regularizers}
\label{sec-examples}

We describe below some popular examples of partly smooth regularizers that are widely used in signal and image processing, statistics and machine learning. We first expose basic building blocks (sparsity, group sparsity, anti-sparsity) and then show how the machinery of partial smoothness enables a powerful calculus to create new priors (using pre- and post-composition, spectral lifting, and positive linear combinations). 

\subsubsection{$\ell^1$ Sparsity}
\label{subsec-l1example}

One of the most popular non-quadratic convex regularization is the $\ell^1$ norm
\eq{
	J(\x) = \norm{x}_1 = \sum_{i=1}^{\N} |\x_i|,
}  
which promotes sparsity. Indeed, it is easy to check that $J$ is partly smooth at $\x$ relative to the subspace 
\eq{
	\Mm_x = T_{\x} = \enscond{ u \in \RR^\N }{ \supp(u) \subseteq \supp(\x) }. 
}
Another equivalent way to interpret this $\ell^1$ prior is that it is the convex enveloppe (restricted to the $\ell^2$-ball) of the $\ell^0$ pseudo-norm~\eqref{eq-lzero-pseudonorm}, in the sense that the $\lun$-unit ball is the convex hull of the restriction of the unit ball of the $\ell^0$-pseudo norm to the $\ldeux$-unit ball.

\myparagraph{Literature review.}
The use of the $\ell^1$ norm as a sparsity-promoting regularizer traces back several decades. An early application was deconvolution in seismology \cite{ClaerboutMuir,santosa1986linear,Taylor}. Rigorous recovery results began to appear in the late 1980's \cite{DonohoStark,DonohoLogan}.
%
In the mid-1990's, $\ell^1$ regularization of least-square problems has been popularized in the signal processing literature under the name basis pursuit~\cite{chen1999atomi} and in the statistics literature under the name Lasso~\cite{tibshirani1996regre}. Since then, the applications and understanding of $\ell^1$ minimization have continued to
increase dramatically.

\subsubsection{$\ell^1-\ell^2$ Group Sparsity}

To better capture the sparsity pattern of natural signals and images, it is useful to structure the sparsity into non-overlapping groups $\Bb$ such that $\bigcup_{b \in \Bb} b = \{1,\ldots,\N\}$. This group structure is enforced by using typically the mixed $\ell^1-\ell^2$ norm
\eql{\label{eq-l1l2-groupsparsity}
	J(\x) = \norm{x}_{1,\Bb} = \sum_{b \in \Bb} \norm{\x_b}, 
} 
where $\x_b = (\x_i)_{i \in b} \in \RR^{|b|}$. Unlike the $\ell^1$ norm, and except the case $|b|=1$ for all $b \in \Bb$, the $\ell^1-\ell^2$ norm is not polyhedral, but is still partly smooth at $\x$ relative to the linear manifold
\eq{
	\Mm_x = T_{\x} = \enscond{ u }{ \supp_\Bb(u) \subseteq \supp_\Bb(\x) }
	\qwhereq
	\supp_\Bb(\x) = \bigcup \enscond{b}{ \x_{b} \neq 0 }.
}

\myparagraph{Literature review.}
The idea of group/block sparsity has been first proposed by~\cite{hall1997numerical,hall1999minimax,cai1999adaptive} for wavelet block shrinkage, i.e. when $\Phi=\Id$. For over-determined regression problems of the form~\eqref{eq:lin-inverse-problem}, it has been introduced by~\cite{bakin1999adaptive,yuan2005model}. Group sparsity has also been extensively used in machine learning in e.g.~\cite{bach2008group} (regression and mutiple kernel learning) and~\cite{obozinski2010joint} (for multi-task learning).The wavelet coefficients of a natural image typical exhibit some group structure, see~\cite{mallat2009a-wav} and references therein on natural image modeling. Indeed, edges and textures induce strong dependencies between coefficients. In audio processing, it has proved useful to structure sparsity in multi-channel data~\cite{gribonval2008atoms}. Group sparsity is also at the heart of the so-called multiple measurements vector (MMV) model, see for instance~\cite{cotter2005sparse,chen2006theoretical}. It is possible to replace the $\ldeux$ norm with more general functionals, such as $\ell^p$ norms for $p > 1$, see for instance~\cite{turlach2005simultaneous,negahban2011simultaneous,vogt2012complete}.

\subsubsection{$\ell^\infty$ Anti-sparsity}

In some cases, the vector to be reconstructed is expected to be flat. Such a prior can be captured using the $\ell^\infty$ norm 
\eq{
	J(x) = \normi{x} = \umax{i \in \ens{1,\dots,n}} \abs{x_i}.
}
It can be readily checked that this regularizer is partly smooth (in fact polyhedral) relative to the subspace
\eq{
 	\Mm_x = T_x = \enscond{u}{u_{I} = \rho x_I \text{ for some } \rho \in \RR}, 
	\qwhereq
	I = \enscond{i}{x_i = \normi{x}} ~.
}

\myparagraph{Literature review.}

The $\linf$ regularization has found applications in computer vision, such as for database image retrieval~\cite{jegou2010improving}. For this application, it is indeed useful to have a compact signature of a signal $x$, ideally with only two values $\pm \normi{x}$ (thus achieving optimal anti-sparsity since $\dim(T_x)=1$ in such a case). An approach proposed in~\cite{jegou2012anti} for realizing this binary quantification is to compute these vectors as solutions of~\eqref{eq-lagrangian} for $J=\norm{\cdot}_\infty$ and a random $\Phi$. A study of this regularization is done in~\cite{fuchs2011spread}, where an homotopy-like algorithm is provided. The use of this $\linf$ regularization is also connected to Kashin's representation~\cite{lyubarskii2010uncertainty}, which is known to be useful in stabilizing the quantization error for instance. Others applications such as wireless network optimization~\cite{studer12signal} also rely on the $\linf$ prior. 


\subsubsection{Synthesis Regularizers}
\label{sec-synth-sparsity}

Sparsity or more general low-complexity regularizations are often used to model coefficients $\al \in \RR^Q$ describing the data $x = D \al$ in a dictionary $D \in \RR^{\N \times Q}$ of $Q$ atoms in $\RR^\N$. 
Given a partly smooth function $J_0 : \RR^Q \rightarrow \RR$, we define the following synthesis-type prior $J : \RR^N \rightarrow \RR$ as the pre-image of $J_0$ under the linear mapping $D$
\eq{
	J(x) = \umin{ \al \in \RR^Q } J_0(\al) \st D\al = x
}
Since $J_0$ is bounded below and convex, $J$ is convex. If $D$ is surjective (as in most cases with redundant dictrionaries), then $J$ is also finite-valued.
The initial optimization~\eqref{eq-lagrangian} can equivalently been solved directly over the coefficients domain to obtain $x^\star = D \al^\star$ where 
\eql{\label{eq-synth-variational}
	\al^\star \in \uArgmin{\al \in \RR^Q} \frac{1}{2\la}\norm{y-\Phi D \al}^2 + J_0(\al)
}
which can be interpreted as a regularized inversion of the operator $\Phi D$ using the prior $J_0$. 

It is possible to study directly the properties of the solutions $\al^\star$ to~\eqref{eq-synth-variational}, which involves directly partial-smoothness of $J_0$. A slightly different question is to understand the behavior of the solutions $x^\star=D\al^\star$ of~\eqref{eq-lagrangian}, which requires to study partial smoothness of $J$ itself. In the case where $D$ is invertible, both problems are completely equivalent.

\myparagraph{Literature review.}
Sparse synthesis regularization using $J_0=\norm{\cdot}_1$ is popular in signal and image processing to model natural signals and images, see for instance~\cite{mallat2009a-wav,starck2010sparse} for a comprehensive account.
The key problem to achieve good performance in these applications is to design a dictionary to capture sparse representations of the data to process. Multiscale dictionaries built from wavelet pyramids are popular to sparsely represent transient signals with isolated singularities and natural images~\cite{mallat1989theory}. The curvelet transform is known to provide non-adaptive near-optimal sparse representation of piecewise smooth images away from smooth edges (so-called cartoon images)~\cite{candes2000curvelets}. Gabor dictionaries (made of localized and translated Fourier atoms) are popular to capture locally stationary oscillating signals for audio processing~\cite{allen1977short}. To cope with richer and diverse contents, researchers have advocated to concatenate several dictionaries to solve difficult problems in signal and image processing, such as component separation or inpainting, see for instance~\cite{elad2005simul}. A line of current active research is to learn and optimize the dictionary from exemplars or even from the available data themselves. We refer to~\cite[Chapter 12]{elad2010sparse} for a recent overview of the relevant literature.

\subsubsection{Analysis Regularizers}
\label{sec-analysis-regularizers}
Analysis-type regularizers (following the terminology introduced in~\cite{elad2007analysis}) are of the form 
\eq{
	J(\x) = J_0(D^* \x) ~,
}
where $D \in \RR^{\N \times Q}$ is a linear operator. Such a prior controls the low complexity (as measured by $J_0$) of the correlations between the columns of $D$ and the signal $x$. If $J_0$ is partly smooth at $z=D^* \x$ for the manifold $\Mm_z^0$, then it is shown in~\cite[Theorem 4.2]{lewis2002active} that $J$ is partly smooth at $\x$ relative to the manifold
\eq{
	\Mm_x = \enscond{ u \in \RR^\N }{ D^* u \in \Mm_z^0 }
}
provided that the following transversality condition holds~\cite[Theorem~6.30(a)]{lee2003smooth}
\eq{
	\Ker(D) \cap \tgtManif{z}{\Mm_z^0}^\perp = \ens{0} \iff \Im(D^*) + \tgtManif{z}{\Mm_z^0} = \RR^N ~.
}

\myparagraph{Literature review.}
A popular example is when $J_0=\norm{\cdot}_1$ and $D^*$ a finite-difference discretization of the derivative of a 1-D signal or a 2-D image. This defines the anisotropic total variation semi-norm, which promotes piecewise constant signals or images~\cite{rudin1992nonlinear}. The 2-D isotropic total variation semi-norm can be interpreted as taking $J_0=\norm{\cdot}_{1,2}$ with blocks of size two. A comprehensive review of total variation regularization can be found in~\cite{chanbolle2009tv}. TV regularization has been extended in several ways to model piecewise polynomial functions, see in particular the Total Generalized Variation prior~\cite{Bredies2010TGV}.


One can also use a wavelet dictionary $D$ which is shift-invariant, such that the corresponding regularization $J$ can be seen as a kind of multi-scale total variation. This is typically the case of the Haar wavelet dictionary~\cite{steidl2004equivalence}. When using higher order wavelets, the corresponding priors favors models $\Mm$ composed of discrete piecewise polynomials. 

The Fused Lasso~\cite{tibshirani2005sparsity} corresponds to $J_0$ being the $\lun$-norm and $D$ is the concatenation of the identity and the adjoint of a finite-difference operator. The corresponding models $\Mm$ are composed of disjoint blocks over which the signals are constant.

Defining a block extracting operator $D^*x = (x_b)_{b \in \Bb}$ allows to re-write the group $\ell^1-\ell^2$ norm~\eqref{eq-l1l2-groupsparsity}, even with overlapping blocks (i.e. $\exists (b,b') \in \Bb^2$ with $b \cap b' \neq \emptyset$), as $J = J_0 \circ D^*$ where $J_0=\norm{\cdot}_{1,2}$ without overlap, see~\cite{jenatton2011structured,2011-eusipco-tree,zhao2009composite,cai2001incorporating}. 
To cope with correlated covariates in linear regression, analysis-type sparsity-enforcing prior were proposed in~\cite{grave2011trace,richard2013inter} using $J_0=\norm{\cdot}_*$ the nuclear norm (as defined in Section~\ref{sec-spectral-priors}).

For unitary $D$, the solutions of~\eqref{eq-lagrangian} with synthesis and analysis regularizations are obviously the same. In the general case (e.g. $D$ overcomplete), however, these two regularizations are different. Some authors have reported results comparing these two priors for the case where $J_0$ is the $\lun$ norm~\cite{elad2007analysis,selesnick2009signal}. A first discussion on the relation and distinction between analysis and synthesis $\lun$-sparse regularizations can be found in~\cite{elad2007analysis}. But only very recently, some theoretical recovery results and algorithmic developments on $\lun$- analysis regularization (so-called cosparse model) have began to be developed, see e.g.~\cite{nam2012cosparse,vaiter2011robust}.

\subsubsection{Spectral Functions}
\label{sec-spectral-priors}

The natural extension of low-complexity priors to matrix-valued data $\x \in \RR^{\N_0 \times \N_0}$ (where $\N=\N_0^2$) is to impose the low-complexity on the singular values of the matrix. We denote $\x = U_{\x} \diag(\La_{\x}) V_{\x}^*$ an SVD decomposition of $\x$, where $\La_{\x} \in \RR_+^{\N_0}$. If $j : \RR^{\N_0} \rightarrow \RR$ is a permutation-invariant closed convex function, then one can consider the function 
\eq{
	J(\x) = j(\La_{\x})
}
which can be shown to be a convex function as well~\cite{LewisMathEig}. When restricted to the linear space of symmetric matrices, $j$ is partly smooth at $\La_{\x}$ for a manifold $m_{\La_{\x}}$, if and only if $J$ is partly smooth at $\x$ relative to the manifold
\eq{
	\Mm_x = \enscond{ U \diag(\La) U^* }{\La \in m_{\La_{\x}}, U \in \Oo_{\N_0}}, 
}
where $\Oo_{\N_0} \subset \RR^{\N_0 \times \N_0}$ is the orthogonal group. The proof of this assertion can be found in~\cite[Theorem~3.19]{daniilidis2013orthogonal}, which builds upon the work of~\cite{DaniilidisMalick14} on manifold smoothness transfer under spectral lifting. This result can be extended to non-symmetric matrices by requiring that $j$ is an absolutely permutation-invariant closed convex function, see~\cite[Theorem~5.3]{daniilidis2013orthogonal}.

\myparagraph{Literature review.}
The most popular spectral prior is obtained for $j=\norm{\cdot}_1$. This defines the nuclear norm, or $1$-Schatten norm, as
\eql{\label{eq-nuclear-norm} 
	J(\x) = \norm{\x}_* = \norm{\La_{\x}}_1 ~.
}
It can be shown that the nuclear norm is the convex hull of the rank function with respect to the spectral norm ball, see~\cite{fazel2002matrix,hiriart2012convex}. It then corresponds to promoting a low-rank prior. Moreover, the nuclear norm can be shown to be partly smooth at $\xx$ relative to the set~\cite[Example~2]{lewis2008alter}
\eq{
	\Mm_x = \enscond{u}{ \rank(u)=\rank(\x) }
}
which is a manifold around $x$.

The nuclear norm has been used in signal and image processing, statistics and machine learning for various applications, including low rank matrix completion~\cite{srebro2004learning,recht2010guaranteed,candes2009exact}, principal component pursuit~\cite{CandesRPCA11}, model reduction~\cite{fazel2001rank}, and phase retrieval~\cite{CandesPhaseLift}. It is also used for some imaging applications, see for instance~\cite{LingalaHDJ11}.

\subsubsection{Mixed Regularizations}

Starting from a collection of convex functions $\{J_\ell\}_{\ell \in \Ll}$, $\Ll=\ens{1,\ldots,L}$, it is possible to design a convex function as 
\eq{
	J_\ell(\x) = \sum_{\ell \in \Ll} \rho_\ell J_{\ell}(\x),
} 
where $\rho_\ell > 0$ are weights. 
If each $J_\ell$ is partly smooth at $\x$ relative to a manifold $\Mm_x^\ell$, then it is shown in~\cite[Corollary 4.8]{lewis2002active} that $J$ is also partly smooth at $\x$ for 
\eq{
	\Mm_x = \bigcap_{\ell \in \Ll} \Mm_x^\ell ~,
}
with the proviso that the manifolds $\Mm_x^{\ell}$ intersect transversally~\cite[Theorem~6.30(b)]{lee2003smooth}, i.e. the sum of their respective tangent spaces $\tgtManif{\x}{\Mm_x^{\ell}}$ spans the whole ambient space $\RR^N$.

\myparagraph{Literature review.}
A popular example is to impose both sparsity and low rank of a matrix, when using $J_1=\norm{\cdot}_1$ and $J_2=\norm{\cdot}_*$, see for instance~\cite{GolbabaeeSparseLowRank,oymak2012simultaneously}.

\subsubsection{Separable Regularization}

Let $\{J_\ell\}_{\ell \in \Ll}$, $\Ll=\ens{1,\ldots,L}$, be a family of convex functions. If $J_\ell$ is partly smooth at $x_\ell$ relative to a manifold $\Mm_{x_\ell}^{\ell}$, then the separable function
\[
J\pa{\{x_\ell\}_{\ell \in \Ll}} = \sum_{\ell \in \Ll} J_{\ell}(\x_\ell)
\]
is partly smooth at $(x_1,\ldots,x_L)$ relative to $\Mm_{x_1}^{1}\times \cdots \times \Mm_{x_L}^{L}$~\cite[Proposition 4.5]{lewis2002active}.

\myparagraph{Literature review.}
One fundamental problem that has attracted a lot of interest in the recent years in data processing involves decomposing an observed object into a linear combination of components/constituents $x_\ell$, $\ell \in \Ll=\ens{1,\ldots,L}$. One instance of such a problem is image decomposition into texture and piece-wise-smooth (cartoon) parts. The corresponding forward model can be cast in the form \eqref{eq:lin-inverse-problem}, where $x_0=\begin{pmatrix}x_1\\x_2\end{pmatrix}$, $x_1$ and $x_2$ are the texture and cartoon components, and $\Phi=[\Id \quad \Id]$. The decomposition is then achieved by solving the variational problem~\eqref{eq-lagrangian}, where $J_1$ is designed to promote the discontinuities in the image, and $J_2$ to favor textures; see \cite{StarckIEEE,Aujol,PeyreFadiliMCA10} and references therein. Another example of decomposition is principal component pursuit, proposed in~\cite{CandesRPCA11}, to decompose a matrix which is the superposition of a low-rank component and a sparse component. In this case $J_1=\norm{\cdot}_1$ and $J_2=\norm{\cdot}_*$.


\section{$\ell^2$ Stability}
\label{sec:intro-robustness}

In this section, we assume that $J$ is a finite-valued convex function, but it is not assumed to be partly smooth.

The observations $y$ are in general contaminated by noise, as described by the forward model~\eqref{eq:lin-inverse-problem}. It is thus important to study the ability of~\eqref{eq-lagrangian} to recover $x_0$ to a good approximation in presence of such a noise $w$, and to assess how the reconstruction error decays as a function of the noise level. In this section, we present a generic result ensuring a so-called ``linear convergence rate'' in terms of $\ell^2$-error between a recovered vector and $x_0$ (see Theorem~\ref{ithm:boundestim}), which encompasses a large body of literature from the inverse problems community.

\subsection{Dual Certificates}

It is intuitively expected that if~\eqref{eq-lagrangian} is good at recovering an approximation of $x_0$ in presence of noise, then~\eqref{eq:reg-noiseless} should be able to identify $x_0$ uniquely when the noise vanishes, i.e. $y=\Phi x_0$. For this to happen, the solution to~\eqref{eq:reg-noiseless} has to satisfy some non-degeneracy condition. To formalize this, we first introduce the notion of dual certificate.

\begin{definition}[Dual certificates]
  	For any vector $x \in \RR^N$, the set of \emph{dual certificates} at $x$ is defined as
  	\eq{
		\Dd(x) = \Im(\Phi^*) \cap \partial J(x) ~. 
	}
\end{definition}
The terminology ``dual certificate'' was introduced in~\cite{candes2009exact}. One can show that the image by $\Phi^*$ of the set of solutions of the Fenchel-Rockafellar dual to~\eqref{eq:reg-noiseless} is precisely $\Dd(x)$. 


It is also worth noting that $x_0$ being a solution of~\eqref{eq:reg-noiseless} for $y=\Phi x_0$ is equivalent to $\Dd(x_0) \neq \emptyset$. Indeed, this is simply a convenient re-writing of the first order optimality condition for~\eqref{eq:reg-noiseless}. 

To ensure stability of the set of minimizers~\eqref{eq-lagrangian} to noise perturbing the observations $\Phi x_0$, one needs to introduce the additional requirement that the dual certificates should be strictly inside the subdifferential of $J$ at $x_0$. This is precisely the non-degeneracy condition mentioned previously.

\begin{definition}[Non-degenerate dual certificates]
  	For any vector $x \in \RR^N$, we define the set of \emph{non-degenerate dual certificates} of $x$
  	\eq{
		\DdN(x) = \Im(\Phi^*) \cap \ri(\partial J(x)) ~. 
	} 
\end{definition}

\subsection{Stability in $\ell^2$ Norm}
\label{sec-linconv-rates}

The following theorem, proved in~\cite{fadili13stable}, establishes a linear convergence rate valid for any regularizer $J$, without any particular assumption beside being a proper closed convex function. In particular, its does not assume partial smoothness of $J$. This generic result encompasses many previous works, as discussed in Section~\ref{sec-litterature-l2}. 

\begin{theorem}\label{ithm:boundestim}
	Assume that 
	\eql{\label{eq-source condition-inj}
		\ker(\Phi) \cap T_{x_0} \cap = \ens{0}
		\qandq
		\DdN(x_0) \neq \emptyset
	}
	and consider the choice $\lambda = c\norm{w}$, for some $c > 0$. Then we have for all minimizers $\xsol$ of~\eqref{eq-lagrangian}
 	\eql{\label{eq-thm-l2-lipsch} 
		\norm{\xsol-x_0}_2 \leq C \norm{w} ~,
  	}
	where $C > 0$ is a constant (see Remark~\ref{rem:l2constant} for details). 
\end{theorem}
In plain words, this bound tells us that the distance of $x_0$ to the set of minimizers of~\eqref{eq-lagrangian} is within a factor of the noise level, which justifies the terminology ``linear convergence rate''. 

\begin{remark}[The role of non-smoothness]
The injectivity of $\Phi$ when restricted to $T_{x_0}$ is intimately related to the fact that $J$ is non-smooth at $x_0$. The higher the degree of non-smoothness, the lower the dimension of the subspace $T_{x_0}$, and hence the more likely the restricted injectivity. If $J$ is smooth around $x_0$ (e.g. quadratic regularizers), however, the restricted injectivity condition cannot be fulfilled, unless $\Phi$ is itself injective. The reason is that $T_{x_0}$ is the whole $\RR^N$ at the smoothness points. For smooth regularizations, it can be shown that the convergence rate is slower than linear, we refer to~\cite{scherzer2009variational} for more details.  
\end{remark}

\begin{remark}[Uniqueness]
One can show that condition~\eqref{eq-source condition-inj} implies that $x_0$ is the unique solution of~\eqref{eq:reg-noiseless} for $y=\Phi x_0$. This condition however does not imply in general that~\eqref{eq-lagrangian} has a unique minimizer for $\la>0$.
\end{remark}

\begin{remark}[Stability constant]\label{rem:l2constant}
The result~\eqref{eq-thm-l2-lipsch} ensures that the mapping $y \mapsto x^\star$ (that might be set-valued) is $C$-Lipschitz-continuous at $y=\Phi x_0$. 
Condition $\DdN(x_0) \neq \emptyset$ is equivalent to the existence of some $\eta \in \DdN(x_0)$. 
The value of $C$ (in fact an upper-bound) can be found in~\cite{fadili13stable}. It depends on $\Phi$, $T_{x_0}$, $c$ and the chosen non-degenerate dual certificate $\eta$. In particular, the constant degrades critically as $\eta$ gets closer to the relative boundary of $\DdN(x_0)$, which reflects the intuition of how far is $\eta$ from being a non-degenerate certificate. 
\end{remark}

\begin{remark}[Source condition]
The condition $\Dd(x_0) \neq \emptyset$ is often called ``source condition'' or ``range condition'' in the literature of inverse problems. We refer to the monograph~\cite{scherzer2009variational} for a general overview of this condition and its implications. It is an abstract condition, which is not easy to check in practice, since exhibiting a valid non-degenerate certificate is not trivial. We give in Section~\ref{sec-linearized-precertif} further insights about this in the context of compressed sensing. Section~\ref{sec-linearized-precertif} describes a particular construction of a good candidate (the so-called linearized pre-certificate) for being such an $\eta \in \DdN(x_0)$, and it is shown to govern stability of the manifold $\Mm_{x_0}$ for partly smooth regularizers. 
\end{remark}

\begin{remark}[Infinite dimension]
It is important to remind that, in its full general form, Theorem~\ref{ithm:boundestim} only holds in finite dimension. 
The constant $C$ indeed may depend on the ambient dimension $N$, in which case the constant can blow-up as the discretization grid of the underlying continuous problem is made finer (i.e. as $N$ grows). 
We detail below some relevant literature where similar results are shown in infinite dimension.
\end{remark}

\subsection{Related Works}
\label{sec-litterature-l2}

\subsubsection{Convergence Rates}

For quadratic regularizations of the form $J = \norm{D^* \cdot}^2$ for some linear operator $D^*$, the $\ldeux$-error decay can be proved to be $O(\sqrt{\norm{w}})$, which is not linear, see~\cite[Chapter~3]{scherzer2009variational} for more details and extensions to infinite dimensional Hilbert spaces. 
For non-smooth priors, in~\cite{burger2004convergence}, the authors show the Bregman distance between $\xsol$ and $x_0$ exhibits a linear convergence rate for both the Lagrangian~\eqref{eq-lagrangian} and the constrained~\eqref{eq-constraint-epsilon} problems under the source condition $\Dd(x_0) \neq 0$. These results hold more generally over infinite dimensional Banach spaces. They have been subsequently generalized to ill-posed non-linear inverse problems by ~\cite{resmerita2005regularization} and~\cite{hofmann2007convergence}.
It is important to observe that in order to prove convergence rates in terms of $\ell^2$-error, as done in~\eqref{eq-thm-l2-lipsch}, it is necessary to strengthen the source condition to its non-degenerate version, i.e. $\DdN(x_0) \neq 0$.

In \cite{lorenz2008convergence}, the authors consider the case where $J$ is a $\ell^p$ norm with $1 \leq p \leq 2$ and establish convergence rates of $\norm{\Phi x_0 - \Phi \xsol}$ in~$O(\norm{w})$ and of $\norm{\xsol - x_0}$ in~$O(\sqrt{\norm{w}})$. 
\cite{grasmair2011necessary} prove Theorem~\ref{ithm:boundestim} for $J=\norm{\cdot}_1$. They show that the non-degeneracy condition is also necessary for linear convergence, and draw some connections with the restricted isometry property (RIP), see below.
Under a condition that bears similarities with~\eqref{eq-source condition-inj}, linear convergence with respect to $J$, i.e. $J(\xsol - x_0) = O(\norm{w})$, is proved in \cite{grasmair2011linear} for positively homogeneous regularizers. This result is equivalent to Theorem~\ref{ithm:boundestim} but only when $J$ is coercive, which precludes many important regularizers, such as for instance analysis-type regularizers including total variation.

\subsubsection{RIP-based Compressed Sensing}

The recovery performance of compressed sensing (i.e. when $\Phi$ is drawn from suitable random ensembles) for $J=\norm{\cdot}_1$ has been widely analyzed under the so-called restricted isometry property (RIP) introduced in~\cite{candes2006robust,candes2006stable,candes2006near}. For any integer $k \geq 0$, the $k$-th order restricted isometry constant of a matrix $\Phi$ is defined as the smallest $\delta_k \geq 0$ such that
\begin{equation*}
 	(1 - \delta_k) \norm{x}^2 \leq \norm{\Phi x}^2 \leq (1 + \delta_k) \norm{x}^2 ,
\end{equation*}
for all vectors $x$ such that $\norm{x}_0 \leq k$. It is shown~\cite{candes2006robust} that if $\delta_{2k} + \delta_{3k} < 1$, then for every vector $x_0$ with $\norm{x_0}_0 \leq k$, there exists a non-degenerate certificate \cite[Lemma 2.2]{candes2005decoding}, see also the discussion in~\cite{grasmair2011necessary}.
In turn, this implies linear convergence rate, and is applied in~\cite{candes2006stable} to show $\ldeux$-stability to noise of compressed sensing.
This was generalized in~\cite{candes2011compressed} to analysis sparsity $J=\normu{D^* \cdot}$, where $D$ is assumed to be a tight frame, structured sparsity in~\cite{candes2011compressed} and matrix completion in~\cite{recht2010guaranteed,candes2011tight} using $J=\norm{\cdot}_*$.
The goal is then to design RIP matrices $\Phi$ with constants such that $\delta_{2k}+\delta_{3k}$ (or a related quantity) is small enough. This is possible if $\Phi$ is drawn from an appropriate random ensemble for some (hopefully optimal) scaling of $(\N,\P,k)$. For instance, if $\Phi$ is drawn from the  standard Gaussian ensemble (i.e. with i.i.d. zero-mean standard Gaussian entries), there exists a constant $C$ such that the RIP constants of $\Phi/\sqrt{\P}$ obey $\delta_{2k}+\delta_{3k}<1$ with overwhelming probability provided that 
\eql{\label{eq-cs-rip-scaling}
	\P \geq C k \log(\N/k) ~,
} 
see for instance~\cite{candes2006near}. This result remains true when the entries of $\Phi$ are drawn independently from a subgaussian distribution. When $\Phi$ is a structured random matrix, e.g. random partial Fourier matrix, the RIP constants of $\Phi/\sqrt{\P}$ can also satisfy the desired bound, but at the expense of polylog terms in the scaling \eqref{eq-cs-rip-scaling}, see~\cite{RauhutFoucart13} for a comprehensive treatment.
Note that in general, computing the RIP constants for a given matrix is an NP-hard problem~\cite{bandeira2013certifying,tillmann2014computational}.

\subsubsection{RIP-less Compressed Sensing}
\label{sec-cs-dual-certif}

RIP-based guarantees are uniform, in the sense that the recovery holds with high probability for \textit{all} sparse signals. There is a recent wave of work in RIP-less analysis of the recovery guarantees for compressed sensing. The claims are non-uniform, meaning that they hold for a fixed signal with high probability on the random matrix $\Phi$. This line of approaches improves on RIP-based bounds providing typically sharper constants. When $\Phi$ is drawn from the Gaussian ensemble, it is proved in~\cite{rudelson2008sparse} for $J=\norm{\cdot}_1$ that if the number of measurements $\P$ obeys $\P \geq C k \log (\N / k)$ for some constant $C>0$, where $k = \norm{x_0}_0$, then condition~\eqref{eq-source condition-inj} holds with high probability on $\Phi$. This result is based on Gordon's comparison principle for Gaussian processes and depends on a summary parameter for convex cones called the Gaussian width. Equivalent lower bounds on the number of measurements for matrix completion from random measurements by minimizing the nuclear norm were provided in~\cite{candes2010power} to ensure that~\eqref{eq-source condition-inj} holds with high probability. This was used to prove $\ldeux$-stable matrix completion in~\cite{candes2010matrix}. 

%

The authors in \cite{chandrasekaran2012convex} have recently showed that the Gaussian-width based approach leads to sharp lower bounds on $\P$ required to solve regularized inverse problems from Gaussian random measurements. For instance, they showed for $J=\norm{\cdot}_1$ that
\eql{\label{eq-cs-gw-scaling}
	\P > 2 k \log(\N/k)
}
guarantees exact recovery from noiseless measurements by solving \eqref{eq:reg-noiseless}. An overhead in the number of measurements is necessary to get linear convergence of the $\ldeux$-error in presence of noise by solving~\eqref{eq-constraint-epsilon} with $\epsilon=\norm{w}$, i.e. $x_0$ is feasible. Their results handle for instance the case of group sparsity~\eqref{eq-l1l2-groupsparsity} and the nuclear norm~\eqref{eq-nuclear-norm}. In the polyhedral case, it can be shown that~\eqref{eq-cs-gw-scaling} implies the existence of a non-degenerate dual certificate, i.e. \eqref{eq-source condition-inj}, with overwhelming probability. The Gaussian width is closely related to another geometric quantity called the statistical dimension in conic integral geometry. The statistical dimension canonically extends the linear dimension to convex cones, and has been proposed in~\cite{LivingOnTheEdge} to deliver reliable predictions about the quantitative aspects of the phase transition for exact noiseless recovery from Gaussian measurements. 

To deal with non-Gaussian matrix measurements (such as for instance partial Fourier matrices), \cite{gross2011recover} introduced the ``golfing scheme'' for noiseless low-rank matrix recovery guarantees using $J=\norm{\cdot}_*$. The golfing scheme is an iterative procedure to construct an (approximate) non-degenerate certificate. This construction is also studied in~\cite{candes2011proba} for noiseless and noisy sparse recovery with $J=\norm{\cdot}_1$. In another chapter of this volume~\cite{TroppChapter14}, the author develops a technique, called the ``bowling scheme'', which is able to deliver bounds on the number of measurements that are similar with the Gaussian width-based bounds for standard Gaussian measurements, but the argument applies to a much wider class of measurement ensembles.




\section{Model Stability}
\label{sec:intro-selection}

In the remainder of this chapter, we assume that $J$ is finite-valued convex and locally partly smooth around $x_0$, as defined in Section~\ref{sec-partial-smoothness}. This means in particular that the prior $J$ promotes locally solution which belong to the manifold $\Mm=\Mm_{x_0}$. In the previous section, we were only concerned with $\ell^2$ stability guarantees and partial smoothness was not necessary then. Owing to the additional structure conveyed by partial smoothness, we will be able to provide guarantees on the identification of the correct $\Mm=\Mm_{x_0}$ by solving~\eqref{eq-lagrangian}, i.e. whether the (unique) solution $x^\star$ of \eqref{eq-lagrangian} satisfies $x^\star \in \Mm$. Such guarantees are of paramount importance for many applications. For instance, consider the case where $\ell^1$ regularization is used to localize some (sparse) sources. Then $x^\star \in \Mm$ means that one perfectly identifies the correct source locations. Another example is that of the nuclear norm for low-rank matrix recovery. The correct model identification implies that $x^\star$ has the correct rank, and consequently that the eigenspaces of $x^\star$ have the correct dimensions and are close to those of $x_0$.

\subsection{Linearized Pre-certificate}
\label{sec-linearized-precertif}

We saw in Section~\ref{sec-linconv-rates} that $\ell^2$-stability of the solutions to~\eqref{eq-lagrangian} is governed by the existence of a non-degenerate dual certificate $p \in \DdN(x_0)$. It turns out that not all dual certificates are equally good for stable model identification, and toward the latter, one actually needs to focus on a particular dual certificate, that we call ``minimal norm'' certificate.

\begin{definition}[Minimal norm certificate]
Assume that $x_0$ is a solution of~\eqref{eq:reg-noiseless}. We define the ``minimal-norm certificate'' as
	\eql{\label{eq-dfn-eta0}
		\eta_0 = \Phi^* \uargmin{\Phi^* p \in \partial J(x_0) } \norm{p} ~.
	}
\end{definition}

A remarkable property, stated in Proposition~\ref{prop-equiv-eta0-etaF} below, is that, as long as one is concerned with checking whether $\eta_0$ is non-degenerate, i.e. $\eta_0 \in \ri(\partial J(x_0))$, one can instead use the vector $\eta_F$ defined below, which can be computed in closed form.

\begin{definition}[Linearized pre-certificate]
	Assume that 
	\eql{\label{eq-inj-etaF}
		\ker(\Phi) \cap T_{x_0} = \{0\}.
	} 
	We define the ``linearized pre-certificate'' as
	\eql{\label{eq-dfn-etaF}
		\eta_F = \Phi^* \uargmin{\Phi^* p \in \aff(\partial J(x_0)) } \norm{p}.
	}	
\end{definition}

\begin{remark}[Well-posedness of the definitions]
	Note that the hypothesis that $x_0$ is a solution of~\eqref{eq:reg-noiseless} is equivalent to saying that $\Dd(x_0)$ is a non-empty convex compact set. Hence in~\eqref{eq-dfn-eta0}, the optimal $p$ is the orthogonal projection of $0$ on a non-empty closed convex set, and thus $\eta_0$ is uniquely defined. Similarly, the hypothesis~\eqref{eq-inj-etaF} implies that the constraint set involved in~\eqref{eq-dfn-etaF} is a non-empty affine space, and thus $\eta_F$ is also uniquely defined.
\end{remark}

\begin{remark}[Certificate vs. pre-certificate]
Note that the only difference between \eqref{eq-dfn-eta0} and \eqref{eq-dfn-etaF} is that the convex constraint set $\partial J(x_0)$ is replaced by a simpler affine constraint. This means that $\eta_F$ does not always qualify as a valid certificate, i.e. $\eta_F \in \partial J(x_0)$, hence the terminology "pre-certificate" used. This condition is actually at the heart of the model identification result exposed in Theorem~\ref{thm-stability}. 
\end{remark}

For now on, let us remark that $\eta_F$ is actually simple to compute, since it amounts to solving a linear system in the least-squares sense.

\begin{proposition}\label{prop-equiv-eta0-etaF}
	Under condition~\eqref{eq-inj-etaF}, one has
	\eql{\label{eq-formula-etaF}
		\eta_F = \Phi^* \Phi_{T_{x_0}}^{+,*} e_{x_0}
		\qwhereq
		e_{x_0} = \proj_{T_{x_0}}(\partial J(x_0)) \in \RR^N.
	} 
\end{proposition}

\begin{remark}[Computating $e_x$]
	The vector $e_x$ appearing in~\eqref{eq-formula-etaF} can be computed in closed form for most of the regularizers discussed in Section~\ref{sec-partial-smoothness}. For instance, for $J=\norm{\cdot}_1$, $e_x=\sign(x)$. For $J=\norm{\cdot}_{1,\Bb}$, it reads $e_x=(e_b)_{b \in \Bb}$, where $e_b=x_b/\norm{x_b}$ if $x_b \neq 0$, and $e_b=0$ otherwise. For $J=\norm{\cdot}_*$ and a SVD decomposition $x=U_{\x} \diag(\La_{\x}) V_{\x}^*$, one has $e_x=U_x V_x^*$.
\end{remark}

The following proposition, whose proof can be found in~\cite{2014-vaiter-ps-consistency}, exhibits a precise relationship between $\eta_0$ and $\eta_F$. In particular, it implies that $\eta_F$ can be used in place of $\eta_0$ to check whether $\eta_0$ is non-degenerate, i.e. $\eta_0 \in \ri(\partial J(x_0))$. 

\begin{proposition}
	Under condition~\eqref{eq-inj-etaF}, one has
	\begin{align*}
		\eta_F \in \ri(\partial J(x_0)) &\qarrq \eta_F = \eta_0, \\
		\eta_0 \in \ri(\partial J(x_0)) &\qarrq \eta_F = \eta_0.
	\end{align*}
\end{proposition}


\subsection{Model Identification}

The following theorem provides a sharp sufficient condition to establish model selection. It is proved in~\cite{2014-vaiter-ps-consistency}. It encompasses as special cases many previous works in the signal processing, statistics ans machine learning literatures, as we discuss in Section~\ref{sec-pw-model-consist}. 

\begin{theorem}\label{thm-stability}
	Let $\J$ be locally partly smooth at $x_0$ relative to $\Mm=\Mm_{x_0}$.
	Assume that
	\eql{\label{eq-condition-model-stability}
		\ker(\Phi) \cap T_{x_0} = \{0\} 
		\qandq
		\eta_F \in \ri(\partial J(x_0)).
	}
	Then there exists $C$ such that if 
	\eql{\label{eq-noise-signal-constr}
		\max(\la, \norm{w}/\la) \leq C, 
	}
	the solution $x^\star$ of~\eqref{eq-lagrangian} from the measurements \eqref{eq:lin-inverse-problem} is unique and satisfies  
  	\eql{\label{eq-model-consistency}
		\xsol \in \Mm \qandq \norm{x_0-\xsol} = O(\max(\la,\norm{w})).
	}
\end{theorem}

\begin{remark}[Linear convergence rate vs. model identification]
Obviously, the assumptions~\eqref{eq-condition-model-stability} of Theorem~\ref{thm-stability} imply those of Theorem~\ref{ithm:boundestim}. They are of course stronger, but imply a stronger result, since uniqueness of $x^\star$ and model identification (i.e. $x^\star \in \Mm$) are not guaranteed by Theorem~\ref{ithm:boundestim} (which does not even need $J$ to be partly smooth). 
A chief advantage of Theorem~\ref{thm-stability} is that its hypotheses can be easily checked and analyzed for a particular operator $\Phi$. Indeed, computing $\eta_F$ only requires solving a linear system, as clearly seen from formula~\eqref{eq-formula-etaF}. 
\end{remark}

\begin{remark}[Minimal signal-to-noise ratio]
Another important distinction between Theorems~\ref{ithm:boundestim} and~\ref{thm-stability} is the second assumption~\eqref{eq-noise-signal-constr}. In plain words, it requires that the noise level is small enough and that the regularization parameter is wisely chosen. Such an assumption is not needed in Theorem~\ref{thm-stability} to ensure linear convergence of the $\ell^2$-error. In fact, this condition is quite natural. To see this, consider for instance the case of sparse recovery where $J=\norm{\cdot}_1$. If the minimal signal-to-noise ratio is low, the noise will clearly dominate the amplitude of the smallest entries, so that one cannot hope to recover the exact support, but it is still possible to achieve a low $\ell^2$-error by forcing those small entries to zero.
\end{remark}

\begin{remark}[Identification of the manifold]
For all the regularizations considered in Section~\ref{sec-examples}, the conclusion of Theorem~\ref{thm-stability} is even stronger as it guarantees that  $\Mm_{\x^\star} = \Mm$. The reason is that for any $\x$ and nearby points $\x'$ with $\x' \in \Mm_{\x}$, one has $\Mm_{\x'}=\Mm_{\x}$.
\end{remark}

\begin{remark}[General loss/data fidelity]\label{rem-generic-loss}
It is possible to extend Theorem~\ref{thm-stability} to account for general loss/data fidelity terms beyond the quadratic one, i.e. $\frac{1}{2}\norm{y-\Phi \x}^2$. More precisely, this result holds true for loss functions of the form $F(\Phi\x,y)$, where $F : \RR^\P \times \RR^\P \rightarrow \RR$ is a $C^2$ strictly convex function in its first argument, $\nabla F$ is $C^1$ in the second argument, with $\nabla F(y,y)=0$, where $\nabla F$ is the gradient with respect to the first variable. In this case, the expression~\eqref{eq-formula-etaF} of $\eta_F$ becomes simply
\eq{
	\eta_F = \Ga ( \proj_T \Ga \proj_T )^+  e_{x_0}	
	\qwhereq
	\choice{
		T = T_{x_0} \\
		\Ga = \Phi^* \partial^2 F(\Phi \x_0,\Phi \x_0) \Phi ~,
	}
}
and where $\partial^2 F$ is the Hessian with respect to the first variable (which is a positive definite operator). We refer to~\cite{2014-vaiter-ps-consistency} for more details.
\end{remark}

\subsection{Sharpness of the Model Identification Criterion}

The following proposition, proved in~\cite{2014-vaiter-ps-consistency}, shows that Theorem~\ref{thm-stability} is in some sense sharp, since the hypothesis $\eta_F \in \ri(\partial J(\x_0))$ (almost) characterizes the stability of $\Mm$.

\begin{proposition}\label{prop-instability}
	We suppose that $\x_0$ is the unique solution of~\eqref{eq:reg-noiseless} for $y=\Phi x_0$ and that
	\eql{\label{eq-instability-cond}
		\ker(\Phi) \cap T_{x_0} = \{0\}, 
		\qandq
		\eta_F \notin \partial J(\x_0).
	}
	Then there exists $C>0$ such that if~\eqref{eq-noise-signal-constr} holds, then any solution $\x^\star$ of~\eqref{eq-lagrangian} for $\la>0$ obeys $\x^\star \notin \Mm$.  
\end{proposition}

In the particular case where $w=0$ (no noise), this result shows that the manifold $\Mm$ is not correctly identified when solving~\eqref{eq-lagrangian} for $y=\Phi x_0$ and for any $\la>0$ small enough. 

\begin{remark}[Critical case]
The only case not covered by neither Theorem~\ref{thm-stability} nor Proposition~\ref{prop-instability} is when $\eta_F \in \rbd(\partial J(\x_0))$, where $\rbd$ stands for the boundary relative to the affine hull. In this case, one cannot conclude, since depending on the noise $w$, one can have either stability or non-stability of $\Mm$. We refer to~\cite{vaiter2011robust} where an example illustrates this situation for the 1-D total variation $J=\norm{D_{\mathrm{DIF}}^* \cdot}_1$, where $D_{\mathrm{DIF}}^*$ is a finite-difference discretization of the 1-D derivative operator.
\end{remark}

\subsection{Probabilistic Model Consistency}
\label{sec-model-consistency-proba}

Theorem~\ref{thm-stability} assumes a deterministic noise $w$, and the operator $\Phi$ is fixed. For applications in statistics and machine learning, it makes sense to rather assume a random model for both $\Phi$ and $w$. The natural question is then to assert that the estimator defined by solving~\eqref{eq-lagrangian} is consistent in the sense that it correctly estimates $x_0$ and possibly the model $\Mm_{x_0}$ as the number of observations $\P \to +\infty$. This requires to handle operators $\Phi$ with an increasing number of rows, and thus to also assess sensitivity of the optimization problem~\eqref{eq-lagrangian} to perturbations of $\Phi$ (and not only to $(w,\la)$ as done previously).  

To be more concrete, in this section, we work under the classical setting where $\N$ an $x_0$ are fixed as the number of observations $\P \to +\infty$. The data $(\Phirow_i,w_i)$ are assumed to be random vectors in $\RR^\N \times \RR$, where $\Phirow_i$ is the $i$-th row of $\Phi$ for $i=1,\ldots,\P$. These vectors are supposed independent and identically distributed (i.i.d.) samples from a joint probability distribution such that $\EE\pa{w_i|\Phirow_i}=0$, finite fourth-order moments, i.e. $\EE\pa{w_i^4} < +\infty$ and $\EE\pa{\norm{\Phirow_i}^4} < +\infty$. Note that in general, $w_i$ and $\Phirow_i$ are not necessarily independent. 
It is possible to consider other distribution models by weakening some of the assumptions and strengthening others, see e.g. \cite{KnightFu2000,zhao2006model,bach2008group}. Let us denote $\Gamma = \EE(\Phirow_i^* \Phirow_i) \in \RR^{\N \times \N}$, where $\Phirow_i$ is any row of $\Phi$. We do not make any assumption on the invertibility of~$\Gamma$. 

In this setting, a natural extension of $\eta_F$ defined by~\eqref{eq-formula-etaF} in the deterministic case is
\eq{
	\tilde\eta_F = \Gamma \Gamma_{T_{x_0}}^{+} e_{x_0} 
}
where $\Gamma_{T_{x_0}}=\proj_{T_{x_0}} \Gamma \proj_{T_{x_0}}$, and we use the fact that $\Gamma_{T_{x_0}}$ is symmetric and $\Im(\Gamma_{T_{x_0}}^{+}) \subset T_{x_0}$. It is also implicitly assumed that $\Ker(\Gamma) \cap T_{x_0} = \ens{0}$ which is the equivalent adaptation of the restricted injectivity condition in~\eqref{eq-condition-model-stability} to this setting. 

To make the discussion clearer, the parameters $(\la=\la_\P,\Phi=\Phi_\P,w=w_\P)$ are now indexed by $\P$. The estimator $\x_\P^\star$ obtained by solving $(\Pp_{\la_\P,y_\P})$ for $y_\P=\Phi_\P x_0 + w_\P$ is said to be consistent for $\x_0$ if, 
\eq{
	\lim_{\P \to +\infty} \Pr\pa{\x_\P^\star ~ \text{is unique}} = 1
} 
and $\x_{\P}^\star \rightarrow \x_0$ in probability. The estimator is said to be model consistent if 
\eq{
	\lim_{\P \to +\infty} \Pr\pa{\x_{\P}^\star \in \Mm} = 1, 
}
where $\Mm=\Mm_{x_0}$ is the manifold associated to $\x_0$.

The following result, whose proof can be found in~\cite{2014-vaiter-ps-consistency}, guarantees model consistency for an appropriate scaling of $\mu_\P$. It generalizes several previous works in the statistical and machine learning literature as we review in Section~\ref{sec-pw-model-consist}. 

\begin{theorem}\label{thm-consistency}
	If 
	\eql{\label{eq-irrepresentable-cond}
		\ker(\Gamma) \cap T_{x_0} = \{0\} 
		\qandq
		\tilde\eta_F \in \ri(\partial J(x_0)), 
	}
	and 
	\eql{\label{eq-lambda-scaling}
		\la_\P = o(\P) 
		\qandq
		\la_\P^{-1} = o(\P^{-1/2}). 
	} 
	Then the estimator $\x_{\P}^\star$ of $\x_0$ is model consistent. 
\end{theorem}

\subsection{Related Works}

\subsubsection{Model Consistency}
\label{sec-pw-model-consist}

Theorem~\ref{thm-stability} is a generalization of a large body of results in the literature. For the Lasso, i.e. $J=\norm{\cdot}_1$, to the best of our knowledge, this result was initially stated in~\cite{fuchs2004on-sp}. In this setting, the result~\eqref{eq-model-consistency} corresponds to the correct identification of the support, i.e. $\supp(\x^\star)=\supp(\x_0)$. 
Condition~\eqref{eq-irrepresentable-cond} for $J=\norm{\cdot}_1$ is known in the statistics literature under the name ``irrepresentable condition'' (generally stated in a non-geometrical form), see e.g. \cite{zhao2006model}. 
\cite{KnightFu2000} have shown estimation consistency for Lasso for fixed $\N$ and $\x_0$ and asymptotic normality of the estimates.
The authors in \cite{zhao2006model} prove Theorem~\ref{thm-consistency} for $J=\norm{\cdot}_1$, though under slightly different assumptions on the covariance and noise distribution.
A similar result is established in \cite{Jia-ElasticNet-Consistency} for the elastic net, i.e. $J=\norm{\cdot}_1 + \rho \norm{\cdot}_2^2$ for $\rho > 0$.
In \cite{bach2008group} and~\cite{bach2008trace}, the author proves Theorem~\ref{thm-consistency} for two special cases, namely the group Lasso and nuclear norm minimization. Note that these previous works assume that the asymptotic covariance $\Corr$ is invertible. We do not impose such an assumption, and only require the weaker restricted injectivity condition $\ker(\Corr) \cap T = \{0\}$. 
In a previous work~\cite{vaiter2011robust}, we have proved an instance of Theorem~\ref{thm-stability} when $J(\x) = \norm{D^* \x}_1$, where $D \in \RR^{\N \times Q}$ is an arbitrary linear operator. This covers as special cases the discrete anisotropic total variation or the fused Lasso. 
This result was further generalized in~\cite{vaiter2013model} when $J$ belongs to the class of partly smooth functions relative to  linear manifolds $\Mm$, i.e. $\Mm=T_{\x}$. Typical instances encompassed in this class are the $\ell^1-\ell^2$ norm, or its analysis version, as well as polyhedral gauges including the $\ell^\infty$ norm. Note that the nuclear norm (and composition of it with linear operators as proposed for instance in~\cite{grave2011trace,richard2013inter}), whose manifold is not linear, does not fit into the framework of~\cite{vaiter2013model}, while it is covered by Theorem~\ref{thm-stability}. 
Lastly, a similar result is proved in~\cite{duval2013spike} for a continuous (infinite dimensional) sparse recovery problem over the space of Radon measures normed by $J$ the total variation of a measure (not to be confused with the total variation of functions). In this continuous setting, an interesting finding is that, when $\eta_0 \in \ri(\partial J(\x_0))$, $\eta_0$ is not equal to $\eta_F$ but to a different certificate (called ``vanishing derivative'' certificate in~\cite{duval2013spike}) that can also be computed by solving a linear system.

\subsubsection{Stronger Criteria for $\lun$}

Many sufficient conditions have been proposed in the literature to ensure that $\eta_F$ is a non-degenerate certificate, and hence to guarantee stable identification of the support (i.e. model). We illustrate this here for $J=\norm{\cdot}_1$, but similar reasoning can be carried out for $\norm{\cdot}_{1,\Bb}$ or $\norm{\cdot}_*$.

The strongest criterion makes use of mutual coherence, first considered in~\cite{donoho2001uncertainty}
\eq{
	\mu(\Phi) = \umax{i \neq j} |\dotp{\phi_i}{\phi_j}| 
}
where each column $\phi_i$ of $\Phi$ are assumed normalized to a unit $\ldeux$ norm. Mutual coherence measures the degree of ill-conditioning of $\Phi$ through the correlation of its columns $(\phi_{i})_{1 \leq i \leq \N}$. Mutual coherence is always lower-bounded by $\sqrt{\tfrac{\N-\P}{\P(\N-1)}}$, and equality holds if and only if $(\phi_{i})_{1 \leq i \leq \N}$ is an equiangular tight frame, see~\cite{strohmer2003grassmannian}. Finer variants based on cumulative coherences have been proposed in~\cite{gribonval2008sparse,borup2008beyond}. To take into account the influence of the support $I=\supp(x_0)$ of the vector $x_0$ to recover, Tropp introduced in~\cite{tropp2006just} the Exact Recovery Condition (ERC), defined as
\eq{
	\text{ERC}(I) = \norm{\Phi_{I^c}^* \Phi_{I}^{+,*}}_{\infty,\infty}
		= \umax{j \notin I} \norm{\Phi_I^+ \phi_j}_1
}
where $\norm{\cdot}_{\infty,\infty}$ is the matrix operator norm induced by the $\ell^\infty$ vector norm, $\Phi_I = (\phi_i)_{i \in I}$ and $I^c$ is the complement of the set $I$. $\Phi_I$ is assumed injective which, in view of Section~\ref{subsec-l1example}, is nothing but a specialization to $\lun$ of the restricted injectivity condition in~\eqref{eq-condition-model-stability}. A weak ERC criterion, which does not involve matrix inversion, is derived in~\cite{DossalSparseSpikes}
\eq{
	\text{wERC}(I) = 
		\frac{       
			\umax{j \in I^c} \sum_{i \in I} \abs{\dotp{\phi_i}{\phi_j}} 
		}{       
			1-\umax{j \in I} \sum_{i \neq j \in I} \abs{\dotp{\phi_i}{\phi_j}}
    	}.
}
Given the structure of the subdifferential of the $\lun$ norm, it is easy to check that
\eq{
	\eta_F \in \ri(\partial J(x_0))
	\quad\Longleftrightarrow\quad	
	\text{IC}(x_0) = \normi{ \Phi_{I^c}^* \Phi_{I}^{+,*} \sign(x_{0,I}) } < 1. 
}
The right hand side in the equivalence is precisely what is called the irrepresentable condition in statistics and machine learning. Clearly, $\text{IC}(x_0)$  involves both the sign vector and the support of $x_0$. The following proposition gives ordered upper bounds of $\text{IC}(x_0)$ in terms of the cruder criteria $\text{ERC}$, $\text{wERC}$ and mutual coherence. A more elaborate discussion of them can be found in~\cite{mallat2009a-wav}. 

\begin{proposition}
Assume that $\Phi_I$ is injective and denote $k=|I|=\norm{x_0}_0$. 
Then,
\eq{
	\text{\upshape IC}(x_0)
	\leq
	\text{\upshape ERC}(I)
	\leq
	\text{\upshape wERC}(I)
	\leq 
	\frac{k \mu(\Phi)}{ 1-(k-1)\mu(\Phi) }.
}
\end{proposition}

\subsubsection{Linearized Pre-certificate for Compressed Sensing Recovery}
\label{sec-cs-modelstab}

Stable support identification has been established in~\cite{wainwright2009sharp,dossal2012sharp} for the Lasso problem when $\Phi$ is drawn from the Gaussian ensemble. These works show that for $k=\norm{x_0}_0$, if 
\eq{
	\P > 2 k \log(\N)
}
then indeed $\eta_F \in \ri(\partial J(x_0))$, and this scaling can be shown to be sharp. This scaling should be compared with~\eqref{eq-cs-gw-scaling} ensuring that there exists a non-degenerate certificate. The gap in the log term indicates that there exists vectors that can be stably recovered by $\ell^1$ minimization in $\ldeux$-error sense, but whose support cannot be stably identified. Equivalently, for these vectors, there exists a non-degenerate certificate but it is not $\eta_F$. 

The pre-certificate $\eta_F$ is also used to ensure exact recovery of a low-rank matrix from incomplete noiseless measurements by minimizing the nuclear norm~\cite{candes2009exact,candes2010power}. This idea is further generalized by~\cite{candes2011simple} for a family of decomposable norms (including in particular $\lun$-$\ldeux$ norm and the nuclear norm), which turns to be a subset of partly smooth regularizers. In these works, lower bounds on the number of random measurements needed for $\eta_F$ to be a non-degenerate certificate are developed. In fact, these measurement lower bounds combined with Theorem~\ref{thm-stability} allow to conclude that matrix completion by solving \eqref{eq-lagrangian} with $J=\norm{\cdot}_*$ identifies the correct rank at high signal-to-noise levels.

\subsubsection{Sensitivity Analysis} 
\label{sec-sensitivity-modelselect}

Sensitivity analysis is a central theme in variational analysis. Comprehensive monographs on the subject are \cite{bonnans2000perturbation,mordukhovich1992sensitivity}.
The function to be analyzed underlying problems~\eqref{eq-lagrangian} and~\eqref{eq:reg-noiseless} is
\eql{\label{eq-parameterizd-functional}
	f(x,\th) = 
	\choice{
		\frac{1}{2\la}\norm{y-\Phi x}^2 + J(x) \qifq \la > 0, \\		
		\iota_{\Hh_y}(x) + J(x) \qifq \la=0, 
	}, 
}
where $\Hh_y = \enscond{y}{\Phi x = y}$ and 
where the parameters are $\th=(\la,y,\Phi)$ for $\la \geq 0$. Theorems~\ref{thm-stability} and~\ref{thm-consistency} can be understood as a sensitivity analysis of the minimizers of $f$ at a point $(\x=\x_0,\th=\th_0=(0,\Phi x_0,\Phi))$.

Classical sensitivity analysis of non-smooth optimization problems seeks conditions to ensure smoothness of the mapping $\th \mapsto \x_\th$ where $x_\th$ is a minimizer of $f(\cdot,\th)$, see for instance~\cite{rockafellar1998var,bonnans2000perturbation}. This is usually guaranteed by the non-degenerate source condition and restricted injectivity condition~\eqref{eq-source condition-inj}, which, as already exposed in Section~\ref{sec-linconv-rates}, ensure linear convergence rate, and hence Lipschitz behavior of this mapping. The analysis proposed by Theorem~\ref{thm-stability} goes one step further, by assessing that $\Mm_{\x_0}$ is a stable manifold (in the sense of~\cite{wright1993ident}), since the minimizer $\x_\th$ is unique and remains in $\Mm_{\x_0}$ for $\th$ close to $\th_0$. Our starting point for establishing Theorem~\ref{thm-stability} is the inspiring work of Lewis~\cite{lewis2002active} who first introduced the notion of partial smoothness and showed that this broad class of functions enjoys a powerful calculus and sensitivity theory.
For convex functions (which is the setting considered in our work), partial smoothness is closely related to $\Uu-\Vv$-decompositions developed in~\cite{Lemarechal-ULagrangian}.
In fact, the behavior of a partly smooth function and of its minimizers (or critical points) depend essentially on its restriction to the manifold, hence offering a powerful framework for sensitivity analysis theory. In particular, critical points of partly smooth functions move stably on the manifold as the function undergoes small perturbations~\cite{lewis2013partial}.  
A important and distinctive feature of Theorem~\ref{thm-stability} is that, partial smoothness of $J$ at $x_0$ relative to $\Mm$ transfers to $f(\cdot,\th)$ for $\la > 0$, but not when $\la=0$ in general. In particular, \cite[Theorem~5.7]{lewis2002active} does not apply to prove our claim.


\section{Sensitivity Analysis and Parameter Selection}
\label{sec:intro-sensitivity}

In this section, we study local variations of the solutions of~\eqref{eq-lagrangian} considered as functions of the observations $y$. In a variational-analytic language, this corresponds to analyzing the sensitivity of the optimal values of~\eqref{eq-lagrangian} to small perturbations of $y$ seen as a parameter. This analysis will have important implications, and we exemplify one of them by constructing unbiased estimators of the quadratic risk, which in turn will allow to have an objectively-guided way to select the optimal value of the regularization parameter~$\la$.

As argued in Section~\ref{sec-sensitivity-modelselect}, assessing the recovery performance by solving~\eqref{eq-lagrangian} for $w$ and $\la$ small amounts to a sensitivity analysis of the minimizers of $f$ in~\eqref{eq-parameterizd-functional} at $(\x=\x_0,\th=\th_0=(0,\Phi x_0,\Phi))$. This section involves again sensitivity analysis of~\eqref{eq-parameterizd-functional} to perturbations of $y$ but for $\la>0$. 
Though we focus our attention on sensitivity to $y$, our arguments extend to any parameters, for instance $\la$ or $\Phi$. 

Similarly to the previous section, we suppose here that $J$ is a finite-valued convex and partly smooth function. For technical reasons, we furthermore assume that the partial smoothness manifold is linear, i.e. $\Mm_x=T_x$. We additionally suppose that the set of all possible models $\Tt = \{T_x\}_{x \in \RR^\N}$ is finite. All these assumptions hold true for the regularizers considered in Section~\ref{sec-examples}, with the notable exception of the nuclear norm, whose manifolds of partial smoothness are non-linear. 

\subsection{Differentiability of Minimizers}
\label{sec:intro-local}

Let us denote $\xsoly$ a minimizer of~\eqref{eq-lagrangian} for a fixed value of $\la > 0$. Our main goal is to study differentiability of $\xsoly$ and find a closed-form formula of the derivative of $\xsoly$ with respect to the observations $y$. Since $\xsoly$ is not necessarily a unique minimizer, such a result means actually that we have to single out one solution $\xsoly$, which hopefully should be a locally smooth function of $y$. However, as $J$ is non-smooth, one cannot hope for such a result to hold for any observation $y \in \RR^\P$. For applications to risk estimation (see Section~\ref{sec-risk-estimation}), it is important to characterize precisely the smallest set $\Hh$ outside of which $\xsoly$ is indeed locally smooth. It turns out that one can actually write down an analytical expression of such a set $\Hh$, containing points where one cannot find locally a smooth parameterization of the minimizers. This motivates our definition of what we coin a ``transition space''.

\begin{definition}[Transition space]
We define the \emph{transition space} $\Hh$ as
\begin{align*}
  \Hh = \bigcup_{\T \in \Tt} \; \bd( \Hh_{\T} ),
\end{align*}
where $\bd(C)$ is the boundary of a set $C$, and
\begin{equation*}
  \Hh_{\T} =
  \enscond{
    y \in \RR^\P
  }{
  	\exists x \in \Tc, \; 
    \la^{-1} \Phi_T^*(\Phi x - y) \in \rbd(\partial J(x))
  }.
\end{equation*}
where $\Tc = \enscond{x \in \RR^\N}{T_x = T}$. 
\end{definition}

The set $\Hh$ contains the observations $y \in \RR^\P$ such that the model subspace $T_{\tilde x(y)}$ associated to a well chosen solution $\tilde x(y)$ of~\eqref{eq-lagrangian} is not stable with respect to small perturbations of $y$. In particular, when $J = \normu{\cdot}$, it can be checked that $\Hh$ is a finite union of hyperplanes and when $J = \norm{\cdot}_{1,2}$ it is a semi-algebraic set (see Definition~\ref{defn-semi-algebraic}). This stability is not only crucial to prove smoothness of $\tilde x(y)$, it is also important to be able to write down an explicit formula for the derivative, as detailed in the following theorem whose proof is given in~\cite{Vaiter-ps-dof}. 

\begin{theorem}\label{ithm-local}
  	Let $y \not\in \Hh$ and $\xxs$ a solution of~\eqref{eq-lagrangian} such that
  	\begin{equation*}\label{eq:introconj-fr}\tag{$\Ii_{\xxs}$}
    	\Ker \Phi_T \cap \Ker \jac^2 J_T(\xxs) = \ens{0}
  	\end{equation*}
  where $\T=\T_{\xxs}$.
  Then, there exists an open neighborhood $\neighb \subset \RR^\N$ of $y$, and a mapping $\solm : \neighb \to \T$ such that
  	\begin{enumerate}
  		\item for every $\bar y \in \neighb$, $\solmB$ is a solution of $(\Pp_{\la,\bar y})$, and $\solm(y) = \xxs$ ; 
  		\item the mapping $\solm$ is $\Calt{1}(\neighb)$ and
    		\begin{equation*}
      			\foralls \bar y \in \neighb, \quad
      			\jac \solm(\bar y) = ( \Phi_T^*\Phi_T + \la \Q_\T(\xxs) )^{-1} \Phi_T .
    		\end{equation*}
  \end{enumerate}
\end{theorem}
Here $\Q_\T$ is the Hessian (second order derivative) of $J$ restricted to $T$. This Hessian is surely well-defined owing to partial smoothness, see~Definition~\ref{dfn-partly-smooth}(\ref{PS-C2}).

\subsection{Semi-algebraic Geometry}

Our goal now is to show that the set $\Hh$ is in some sense ``small'' (in particular to show that it has zero Lebesgue measure), which will entail differentiability of $y \mapsto \xxs$ Lebesgue almost everywhere. For this, additional geometrical structure on $J$ is needed. Such a rich class of functions is provided by the notion of a semi-algebraic subset of $\RR^\NN$ to be defined shortly. Semi-algebraic sets and functions have been broadly applied to various areas of optimization. The wide applicability of semi-algebraic functions follows largely from their stability under many mathematical operations. In particular, the celebrated Tarski-Seidenberg theorem states, loosely, that the projection of a semi-algebraic set is semi-algebraic. These stability properties are crucial to obtain the following result, proved in~\cite{Vaiter-ps-dof}. 

\begin{definition}[Semi-algebraic set and function]\label{defn-semi-algebraic}
A set $E$ is semi-algebraic if it is a finite union of sets defined by polynomial equations and (possibly strict) inequalities. A function $f : E \rightarrow F$ is semi-algebraic if $E$ and its graph $\enscond{(u,f(u))}{u \in E}$ are semi-algebraic sets. 
\end{definition}

\begin{remark}[From semi-algebraic to o-minimal geometry]
The class of semi-alge\-braic functions is large, and subsumes, for instance, all the regularizers $J$ described in Section~\ref{sec-examples}. The qualitative properties of semi-algebraic functions are shared by a much bigger class called functions definable in an o-minimal structure over $\RR$, or simply definable functions. O-minimal structures over $\RR$ correspond in some sense to an axiomatization of some of the prominent geometrical properties of semi-algebraic geometry \cite{coste1999omin} and particularly of the stability under projection. For example, the function $J(x)=\sum_i |x_i|^s$, for an arbitrary $s \geq 0$, is semi-algebraic only for rational $s \in \QQ$, while it is always definable in an o-minimal structure \cite{vandendries1996geom}. Due to the variety of regularizations $J$ that can be formulated within the framework of o-minimal structures, all our results stated in this section apply to definable functions, see~\cite{Vaiter-ps-dof} for a detailed treatment.
\end{remark}

Semi-algebraic functions are stable for instance under (sub)differentiation and projection. These stability properties are crucial to obtain the following result, proved in~\cite{Vaiter-ps-dof}. 

\begin{proposition}
\label{prop:transsemialgebraic}
If $J$ is semi-algebraic, the transition space $\Hh$ is semi-algebraic and has zero Lebesgue measure.
\end{proposition}

\subsection{Unbiased Risk Estimation}
\label{sec-risk-estimation}

A problem of fundamental practical importance is to automatically adjust the parameter $\la$ to reach the best recovery performance when solving~\eqref{eq-lagrangian}. Parameter selection is a central theme in statistics, and is intimately related to the question of model selection, as introduced in Section~\ref{sec-model-selection}.

We then adopt a statistical framework in which the observation model~\eqref{eq:lin-inverse-problem} becomes
\eql{\label{eq-fwd-random}
	Y = \Phi x_0 + W
}  
where $W$ is random noise having an everywhere strictly positive probability density function, assumed to be known. Though the forthcoming results can be stated for a large family of distributions, for the sake of concreteness, we only consider the white Gaussian model where $W \sim \Nn(0,\sigma^2 \Id_{\P \times \P})$, with known variance $\sigma^2$.

Under the observation model~\eqref{eq-fwd-random}, the ideal choice of $\la$ should be the one which minimizes the quadratic estimation risk $\EE_W(\norm{\xsol(Y)-x_0}^2)$. This is obviously not realistic as $x_0$ is not available, and in practice, only one realization of $Y$ is observed. To overcome these obstacles, the traditional approach is to replace the quadratic risk with an some estimator that solely depends on $Y$. The risk estimator is also expected to enjoy nice statistical properties among which unbiasedness is highly desirable.

However, it can be shown, see e.g.~\cite[Section~IV]{eldar2009gsure}, that the quadratic risk $\EE_W(\norm{\xsol(Y)-x_0}^2)$ cannot be reliably estimated on $\ker(\Phi)$. Nonetheless, we may still obtain a reliable assessment of the part that lies in $\Im(\Phi^*)=\ker(\Phi)^\perp$ or any linear image of it. For instance, the most straightforward surrogate of the above risk is the so-called prediction risk $\EE_W(\norm{\mu(Y)-\mu_0}^2)$, where 
\eq{
	\mu_0=\Phi x_0
	\qandq
	\msol(y) = \Phi \xsoly,
}
where $\xsoly$ is any solution of~\eqref{eq-lagrangian}.
One can easily show that $\msol(y) \in \RR^\P$ is well-defined as a single-valued mapping and thus does not depend on the particular choice of $\xsoly$, see~\cite{Vaiter-ps-dof}. Consequently, Theorem~\ref{ithm-local} shows that $y \mapsto \msol(y)$ is a $C^1$ mapping on $\RR^\P \setminus \Hh$. 

\subsection{Degrees of Freedom}
\label{sec-dof}

The degrees of freedom (DOF) quantifies the model ``complexity'' of a statistical modeling procedure \cite{efron1986biased}. It is at the heart of several risk estimation procedures. Therefore, in order to design estimators of the prediction risk, an important step is to get an estimator of the corresponding DOF.

\begin{definition}[Empirical DOF]
	Suppose that $y \mapsto \msol(y)$ is differentiable Lebesgue almost everywhere, as is the case when it is Lipschitz-continuous (Rademacher's theorem). The empirical number of degrees of freedom is defined as
	\eq{
		\DOF(y) = \diverg(\msol)(y) = \tr(\jac \msol(y)), 
	}
	where the derivative is to be understood in the weak sense, i.e. to hold Lebesgue almost everywhere (a.e.).
\end{definition}

An instructive example to get the gist of this formula is the case where $\msol$ is the orthogonal projection onto some linear subspace $V$. We then get easily that $\DOF(y) = \dim(V)$, which is in agreement with the intuitive notion of the number of DOF.

The following result delivers the closed-form expression of $\DOF(y)$, valid on a full Lebesgue measure set, for $\msol(y)=\Phi\xsol(y)$ and $\xsol(y)$ an appropriate solution of~\eqref{eq-lagrangian}. At this stage, it is important to realize that the main difficulty does not lie in showing almost everywhere differentiability of $\msoly$; this mapping is in fact Lipschitz-continuous by classical arguments of sensitivity analysis applied to~\eqref{eq-lagrangian}. Rather, it is the existence of such a formula and its validity Lebesgue a. e. that requires more subtle arguments obtained owing to partial smoothness of $J$. For this, we need also to rule out the points $y$ where~\eqref{eq:introconj-fr} does not hold. This is the rationale behind the following set.

\begin{definition}[Non-injectivity set]
We define the \emph{Non-injectivity set} $\Gg$ as
\begin{align*}
  \Gg = \enscond{y \notin \Hh}{\text{\eqref{eq:introconj-fr} does not hold for any minimizer $\xxs$ of~\eqref{eq-lagrangian}}} ~.
\end{align*}
\end{definition}


\begin{theorem}\label{ithm-dof}
  For every $y \notin \Hh \cup \Gg$, there is $\xxs$ such that~\eqref{eq-lagrangian} holds and
  \eql{\label{eq-analytic-dof}
    	\DOF(y) = \tr(\Delta_{x^\star}(y))
    	\qwhereq
     	\Delta_{x^\star}(y) = \XX_\T \circ ( \transp{\XX_\T} \XX_\T + \la\Q_\T(\xxs) )^{-1} \circ \transp{\XX_\T}, 
  }
  were $T=T_{x^\star}$. 
\end{theorem}

\begin{remark}[Non-injectivity set]
It turns out that $\Gg$ is in fact empty for many regularizers. This is typically the case for $J=\norm{\cdot}_1$~\cite{dossal2013degrees}, $J=\norm{D^*\cdot}_1$~\cite{vaiter2013local}, and the underlying reasoning can be more generally extended to polyhedral regularizers. The same result was also shown for $J=\norm{\cdot}_{1,2}$ in~\cite{vaiter2013degrees}. More precisely, in all these works, it was shown that for each $y \notin \Hh$, there exists a solution $\xsol$ of~\eqref{eq-lagrangian} that fulfills \eqref{eq:introconj-fr}. The proof is moreover constructive allowing to build such a solution starting from any other one.
\end{remark}

\subsection{Stein Unbiased Risk Estimator (SURE)}
\label{sec-SURE}

We now have all necessary ingredients at hand to design an estimator of the prediction risk.

\begin{definition}
	Suppose that $y \mapsto \msol(y)$ is differentiable Lebesgue almost everywhere, as is the case when it is Lipschitz-continuous. The SURE is defined as
	\eql{\label{eq-sure-def}
		\SURE(y) = \norm{y - \msol(y)}^2 + 2 \sigma^2 \DOF(y) - \P \sigma^2.
	}
\end{definition}

In this definition, we have anticipated on unbiasedness of this estimator. In fact, this turns out to be a fundamental property owing to the celebrated lemma of Stein~\cite{stein1981estimation}, which indeed asserts that the SURE \eqref{eq-sure-def} is an unbiased estimator of the prediction risk. Therefore, putting together Theorem~\ref{ithm-dof}, Proposition~\ref{prop:transsemialgebraic} and Stein's lemma, we get the following.

\begin{theorem}
	Suppose that $J$ is semi-algebraic and $\Gg$ is of zero Lebesgue measure. Then,
	\eq{
		\EE_W( \SURE(Y) ) = \EE_W(\norm{\mu(Y)-\mu_0}^2)
	}
	where \eqref{eq-analytic-dof} is plugged into \eqref{eq-sure-def}, and $\mu(Y)=\Phi\xsol(Y)$.
\end{theorem}


\begin{remark}[Parameter selection]
A practical usefuleness of the SURE is its ability to provide an objectively guided way to select a good $\la$ from a single observation $y$ by minimizing $\SURE(y)$. While unbiasedeness of the SURE is guaranteed, it is hard to control its variance and hence its consistency. This is an open problem in general, and thus little can be said about the actual theoretical efficiency of such an empirical parameter selection method. It works however remarkably well in practice, see the discussion in Section~\ref{sec-SURE-numerics} and references therein.
\end{remark}

\begin{remark}[Projection risk]
The SURE can be extended to unbiasedly estimate other risks that the prediction one. For instance, as argued in Section~\ref{sec-risk-estimation}, one can estimate the so-called projection risk defined as $\EE_W( \|\proj_{\ker(\Phi)^\bot}(\xsol(Y)-x_0) \|^2)$. This is obviously better that the prediction risk as a surrogate for the estimation risk. 
\end{remark}

\subsection{Related Works}

\subsubsection{Sensitivity Analysis}

In Section~\ref{sec-sensitivity-modelselect}, we reviewed the relevant literature pertaining to sensitivity analysis for partly smooth functions, which is obviously very connected to Theorem~\ref{ithm-local}. See also \cite{bolte2011generic} for the case of linear optimization over a convex semi-algebraic partly smooth feasible set, where the authors prove a sensitivity result with a zero-measure transition space. A distinctive feature of our analysis toward proving unbiasedness of the SURE is the need to ensure that sensitivity analysis can be carried out on a full Lebesgue measure set. In particular, it necessitates local stability of the manifold $\Mm_{\xxs}$ associated to an appropriate solution $\xxs$, and this has to hold Lebesgue almost everywhere. Thus the combination of partial smoothness and semi-algebraicity is the key.

\subsubsection{Risk Estimators}

In this section, we put emphasis on the SURE as an unbiased estimator of the prediction risk. There are other alternatives in the literature which similarly rely on estimator of the DOF. One can think for instance of the generalized cross-validation (GCV)~\cite{golub1979generalized}. Thus our results apply equally well to such risk estimators.
Extensions of the SURE to independent variables from a continuous exponential family are considered in \cite{hudson1978nie}. \cite{eldar2009gsure} generalizes the SURE principle to continuous multivariate exponential families, see also~\cite{pesquet-deconv,vaiter2013local} for the multivariate Gaussian case. The results described here can be extended to these setting as well, see~\cite{Vaiter-ps-dof}.

\subsubsection{Applications of SURE in Statistics and Imaging}

Applications of SURE emerged for choosing the parameters of linear estimators such ridge regression or smoothing splines~\cite{li1985sure}.
After its introduction in the wavelet community through the SURE-Shrink estimator \cite{donoho1995adapting}, it has been extensively used for various image restoration problems, e.g.~%
with sparse regularization \cite{blu2007surelet,vonesch2008sure,ramani2008montecarlosure,chaux2008nonlinear,pesquet-deconv,cai2009data,luisier2010sure,ramani2012regularization,ramani2012iterative} or with non-local means~\cite{vandeville2009sure,duval2011abv,dds2011nlmsap,vandeville2011non}.

\subsubsection{Closed-form Expressions for SURE}

For the Lasso problem, i.e. $J=\norm{\cdot}_1$, the divergence formula~\eqref{eq-analytic-dof} reads
\begin{equation*}
  \DOF(y) = \abs{\supp(\xxs)},
\end{equation*}
where $\xxs$ is a solution of~\eqref{eq-lagrangian} such that~\eqref{eq:introconj-fr} holds, i.e. $\Phi_{\supp(\xxs)}$ has full rank.
This result is proved in~\cite{zou2007degrees} for injective $\Phi$ and in~\cite{dossal2013degrees} for arbitrary $\Phi$. This result is extended to analysis $\lun$-sparsity, i.e. $J=\norm{D^* \cdot}_1$, in~\cite{tibshirani2012degrees,vaiter2013local}. A formula for the DOF in the case where $\xsoly$ is the orthogonal projection onto a partly smooth convex set $C$ is proved in~\cite{kato2009degrees}. This work extends that of~\cite{meyer2000degrees} which treats the case where $C$ is a convex polyhedral cone. These two works allow one to compute the degrees of freedom of estimators defined by solving~\eqref{eq-constraint-gamma} in the case where $\Phi$ is injective. \cite{hansen2014dof} studied the DOF of the metric projection onto a closed set (non-necessarily convex), and gave a precise representation of the bias when the projection is not sufficiently differentiable.

A formula of an estimate of the DOF for the group Lasso, i.e. $J=\norm{\cdot}_{1,2}$ when $\Phi$ is orthogonal within each group was conjectured in~\cite{yuan2005model}. An estimate is also given by~\cite{solo2010threshold} using heuristic derivations that are valid only when $\Phi$ is injective, though its unbiasedness is not proved. \cite{vaiter-icml-workshops} derived an estimator of the DOF of the group Lasso and proved its unbiasedness when $\Phi$ is injective. Closed-form expression of the DOF estimate for denoising with the nuclear norm, i.e. $\Phi=\Id$ and $J=\norm{\cdot}_*$, were concurrently provided in~\cite{DeledalleSVTSURE12,CandesSVTSURE12}.

\subsubsection{Numerical Methods for SURE}
\label{sec-SURE-numerics}

Deriving the closed-form expression of the DOF is in general challenging and has to be addressed on a case by case basis. The implementation of the divergence formula such as~\eqref{eq-analytic-dof} can be computationally expensive in high dimension. But since only the trace of the Jacobian is needed, it is possible to speed up these computations through Monte-Carlo sampling, but at the price of mild approximations.
If the Jacobian is not known in closed-form or prohibitive to compute, one may appeal to finite-difference approximations along Monte Carlo sampled directions~\cite{ye1998measuring,shen2002adaptive}, see~\cite{girard1989fast,ramani2008montecarlosure} for applications to imaging problems.

In practice, the analytical formula~\eqref{eq-analytic-dof} might be subject to serious numerical instabilities, and thus cannot always be applied safely when the solution $\xxs$ is only known approximately. Think for instance of the case where $\xxs$ is approximated by an an iterate computed after finitely-many iterations of an algorithm as detailed in Section~\ref{sec-algorithms}. A better practice is then to directly compute the DOF, hence the SURE, recursively from the iterates themeselves, as proposed by~\cite{vonesch2008sure,giryes-proj-gsure,deledalle2012proximalncmip}.


\section{Proximal Splitting for Structured Optimization}
\label{sec-algorithms}

Though problems~\eqref{eq-lagrangian},~\eqref{eq:reg-noiseless},~\eqref{eq-constraint-epsilon} or~\eqref{eq-constraint-gamma} are non-smooth, they enjoy enough structure to be solved by efficient algorithms. The type of algorithm to be used depends in particular on the properties of $J$. We first briefly mention some popular non-smooth optimization schemes in Section~\ref{sec-cvx-optim}, and focus our attention on proximal splitting schemes afterwards.

\subsection{Convex Optimization for Regularized Inverse Problems}
\label{sec-cvx-optim}

\subsubsection{(Sub)-gradient Descent} 

Consider for example problem~\eqref{eq-lagrangian}. This is a convex composite optimization problem where one of the functions is smooth with a Lipschitz-continuous gradient. If $J$ were smooth enough, then a simple gradient (or possibly (quasi-)Newton) descent method could be used. However, as detailed in Section~\ref{sec-partial-smoothness}, low-complexity regularizers $J$ are intended to be non-smooth in order to promote models $\Mm$ of low intrinsic dimension, and $J$ is precisely non-smooth transverse to $\Mm$. One can think of replacing gradients by subgradients (elements of the subdifferential), since $J$ is assumed finite-valued (hence closed) convex, which are bounded. This results in a subgradient descent algorithm which is guaranteed to converge but under stringent assumptions on the descent step-sizes, which in turn makes their global convergence rate quite slow, see~\cite{NesterovSmooth}. 


\subsubsection{Interior Point Methods} 

Clearly, the key to getting efficient algorithms is to exploit the structure of the optimization problems at hand while handling non-smoothness properly. For a large class of regularizers $J$, such as those introduced in Section~\ref{sec-examples}, the corresponding optimization problems can be cast as conic programs. The cone constraint can be enforced using a self-concordant barrier function, and the optimization problem can hence be solved using interior point methods, as pioneered by~\cite{nesterov1994interior}, see also the monograph~\cite{boyd2004convex}. This class of methods enjoys fast convergence rate. Each iteration however is typically quite costly and can become prohibitive as the dimension increases.

\subsubsection{Conditional Gradient} 

This algorithm is historically one of the first method for smooth constrained convex optimization (a typical example being~\eqref{eq-constraint-gamma}), and was extensively studied in the 70's. It is also known as Frank-Wolfe algorithm, since it was introduced by~\cite{FrankWolfe} for quadratic programming and extended in~\cite{DunnHarshbarger}.
The conditional gradient algorithm is premised on being able to easily solve (at each iteration) linear optimization problems over the feasible region of interest. This is in contrast to other first-order methods, such as forward-backward splitting and its variants (see Section~\ref{sec-fb}), which are premised on being able to easily solve (at each iteration) a projection problem. Moreover, in many applications the solutions to the linear optimization subproblem are highly structured and exhibit particular sparsity and/or low-rank properties. These properties have renewed interest in the conditional gradient method to solve sparse recovery ($\ell_1$ and total variation), low-rank matrix recovery (nuclear norm minimization), anti-sparsity recovery, and various other problems in signal processing and machine learning; see e.g.~\cite{Clarkson08,Sulovsky10,Shalev11,DudikHM12,Juditsky14}.


\subsubsection{Homotopy/Path-following} 

Homotopy and path-following-type methods have been introduced in the case of $\lun$-minimization to solve~\eqref{eq-lagrangian} by~\cite{osborne2000new}. They were then adapted to analysis $\lun$, i.e. $J=\normu{D^*\cdot}$, in~\cite{tibshirani2011solution}, and $\linf$ regularization, $\normi{\cdot}$, in~\cite{fuchs2011spread}. One can in fact show that these methods can be applied to any polyhedral regularization (see~\cite{vaiter13polyhedral}), because these methods only rely on the crucial fact that the solution path $\lambda \mapsto \xxs_\la$, where $\xxs_\la$ is a solution of~\eqref{eq-lagrangian}, is piecewise affine.
The LARS algorithm~\cite{efron2004least} is an accelerated version of homotopy which computes an approximate homotopy path for $J=\norm{\cdot}_1$ along which the support increases monotonically along the course of iterations.
In the noiseless compressed sensing case, with $\Phi$ drawn from the Gaussian ensemble, it is shown in~\cite{donoho2008fast} that if $x_0$ is $k$-sparse with $\P > 2 k \log(\N)$, the homotopy method reaches $x_0$ in only $k$ iterations. This $k$-solution property was empirically observed for other random matrix ensembles, but at different thresholds for $\P$. In~\cite{MairalYu12}, the authors proved that in the worst case, the number of segments in the solution path is exponential in the number of variables, and thus the homotopy method can then take as many iterations to converge. 

As for interior points, the cost per iteration of homotopy-like methods, without particular ad hoc optimization, scales badly with the dimension, thus preventing them to be used for large-scale problems such as those encountered in imaging. This class of solvers is thus a wise choice for problems of medium size, and when high accuracy (or even exact computation up to machine precision for the homotopy algorithm) is needed. Extensions of these homotopy methods can deal with progressive changes in the operator $\Phi$ or the observations $y$, and are thus efficient for these settings, see~\cite{AsifHomotopy}.

\subsubsection{Approximate Message Passing} 

In the last five years, ideas from graphical models and message passing and approximate message passing algorithms have been proposed to solve large-scale problems of the form~\eqref{eq-lagrangian} for various regularizers $J$, in particular $\lun$, $\lun-\ldeux$ and the nuclear norm. A comprehensive review is given in~\cite{MontanariReview12}. However, rigorous convergence results have been proved so far only in the case in which $\Phi$ is standard Gaussian, though numerical results show that the same behavior should apply for broader random matrix ensembles.

\subsection{Proximal Splitting Algorithms}

Proximal splitting methods are first-order iterative algorithms that are tailored to solve structured non-smooth (essentially convex) optimization problems. The first operator splitting method has been developed from the 70's. Since then, the class of splitting methods have been regularly enriched with increasingly sophisticated algorithms as the structure of problems to handle becomes more complex. 

%
%
To make our discussion more concrete, consider the general problem of minimizing the proper closed convex function 
\eq{\label{eq:minpbgen}
	f = h + \sum_{k=1}^K g_k \circ A_k
}
where $h : \RR^N \rightarrow \RR$ is convex and smooth, the $A_k : \RR^N \rightarrow \RR^{N_k}$ are linear operators and $g_k : \RR^{N_k} \rightarrow \RR$ are proper closed convex functions for which the so-called proximity operator (to be defined shortly) can be computed easily (typically in closed form). We call such a function $g_k$ ``simple''.

\begin{definition}
The proximity operator of a proper closed convex function $g$ is defined as, for $\ga > 0$
\eq{
	\Prox_{\ga g}(\x) = \uargmin{u \in \RR^\N} \frac{1}{2}\norm{x-u}^2 + \ga g(u).
}
\end{definition}
The proximal operator generalizes the notion of orthogonal projection onto a non-empty closed convex set $C$ that one recovers by taking $g = \iota_C$. 

Proximal splitting algorithms may evaluate (possibly approximately) the individual operators (e.g. gradient of $h$), the proximity operators of the $g_k$'s, the linear operators $A_k$, all separately at various points in the course of iteration, but never those of sums of functions nor composition by a linear operator. Therefore, each iteration is cheap to compute for large-scale problems. They also enjoy rigorous convergence guarantees, stability to errors, with possibly quantified convergence rates and iteration complexity bounds on various quantities. This justifies their popularity in contemporary signal and image processing or machine learning, despite that their convergence is either sublinear or at best linear.

It is beyond the scope of this Chapter to describe thoroughly the huge literature on proximal spliting schemes, as it is a large and extremely active research field in optimization theory. Good resources and reviews on the subject are~\cite{beck2009gradient,bauschke2011convex,combettes2011proximal,parikh2013proximal}. We instead give a brief classification of the most popular algorithms according to the class of structured objective functions they are able to handle:
\begin{itemize}
	\item \textit{Forward-Backward (FB)} algorithm~\cite{Mercier79,passty1979ergodic,combettes2005signal}. It is designed to minimize~\eqref{eq:minpbgen} when $h$ has a Lipschitz-continuous gradient, $K=1$, $A_1=\Id$, and $g_1$ is simple. There are accelerated (optimal) variants of FB, such as the popular Nesterov~\cite{Nesterov07} or Fista~\cite{beck2009fista}, but the convergence of the iterates is not longer guaranteed for these schemes. FB and its variants are good candidates to solve~\eqref{eq-lagrangian}. We will further elaborate on FB in Section~\ref{sec-fb}.
	\item \textit{Douglas-Rachford (DR)} algorithm~\cite{douglas1956numerical,lions1979splitting}. It is designed to minimize~\eqref{eq:minpbgen} for $h=0$, $K=2$, $A_k=\Id$ and $g_k$ is simple for $k=1,2$. It can be easily extended to the case of $K > 2$ by either lifting to a product space, see e.g.~\cite{CombettesPesquet08}, or through projective splitting~\cite{EcksteinSvaiter09}. DR can be used to solve~\eqref{eq:reg-noiseless},~\eqref{eq-constraint-epsilon} or~\eqref{eq-constraint-gamma} for certain operators $\Phi$.	
	\item \textit{Generalized Forward-Backward (GFB)} algorithm~\cite{gfb2011}. It can handle the case of an arbitrary $K$ with $A_k=\Id$, $g_k$ simple and $h$ has a Lipschitz-continuous gradient. It can be interpreted as hybridization of FB scheme and the DR scheme on a product space. 
	\item \textit{Alternate Direction Method of Multipliers (ADMM)} algorithm~\cite{fortin2000augmented,gabay1983chapter,gabay1976dual,glowinski1989augmented}. It is adapted to minimize~\eqref{eq:minpbgen} for $h=0$, $K=2$ with $A_1=\Id$ and $A_2$ is injective. It can be shown~\cite{gabay1976dual,eckstein1992douglas} that ADMM is equivalent to DR applied to the Fenchel-Rockafellar dual problem $\min_u g_1^* \circ -A_2^*(u) + g_2^*(u)$, where $g_k^*$ is the Legendre-Fenchel conjugate of $g_k$. While DR applies when $g_1$ and $g_2 \circ A_2$ are simple, ADMM is a better alternative whereas both $g_1 \circ -A_2^*$ and $g_2^*$ are simple. Extension to the case $K > 2$ was proposed for instance in~\cite{Eckstein94}.
	\item \textit{Dykstra} algorithm~\cite{Dykstra83}. It is able to solve the case where $h(x)=\norm{x-y}^2$, $A_k=\Id$ and the $g_k$ are simple functions. It was initially introduced by~\cite{Dykstra83} in the case where the $g_k$ are indicator functions of closed convex sets, and is generalized in~\cite{BauschkeCombettes-Dykstra} to arbitrary convex functions. It is also extended in~\cite{CensorReich-Dykstra,BauschkeLewis-dykstras} to the case where $h$ is a Bregman divergence. 
	\item \textit{Primal-Dual schemes}. Recently, primal-dual splitting algorithms have been proposed to minimize~\eqref{eq:minpbgen} in its full generality, and even more complex objectives, see for instance~\cite{chen1994proximal,tseng1997alternating,solodov2004class,briceno2011monotone,chambolle2011first,combettes2012primal,condat2012primal,vu2011splitting}. Primal-dual schemes can be used to solve~\eqref{eq-lagrangian},~\eqref{eq:reg-noiseless},~\eqref{eq-constraint-epsilon} or~\eqref{eq-constraint-gamma}.
\end{itemize}

\subsection{Finite Model Identification with Forward Backward}
\label{sec-fb}

\newcommand{\iter}[1]{{#1}^{(n)}}
\newcommand{\siter}[1]{{#1}_{n}}
\newcommand{\iit}[1]{{#1}^{(n+1)}}

The FB algorithm is a good candidate to solve~\eqref{eq-lagrangian} when $J$ is simple. Starting from some $x^{(0)} \in \RR^\N$, the FB iteration applied to ~\eqref{eq-lagrangian} reads
\eq{
	\iit{x} = \text{Prox}_{\siter{\tau}\la J}\pa{ \iter{x} + \siter{\tau} \Phi^* ( y - \Phi \iter{x} ) },  
}
where the step-size sequence should satisfy $0 < \underline{\tau} \leq \siter{\tau} \leq \overline{\tau} < 2/\norm{\Phi}^2$ to ensure convergence of of the sequence $\iter{x}$ to a minimizer of~\eqref{eq-lagrangian}.

In fact, owing to partial smoothness of $J$, much more can be said about the iterates of the FB algorithm. More precisely, after a finite number of iterations, Forward-Backward algorithm correctly identifies the manifold $\Mm$. This is made formal in the following theorem whose proof can be found in~\cite{Liang14}. 

\begin{theorem}\label{thm-fb}
Under the assumptions of Theorem~\ref{thm-stability}, $\iter{x} \in \Mm$ for $n$ large enough.
\end{theorem}

This result sheds some light on the convergence behavior of this algorithm in the favorable case where condition~\eqref{eq-inj-etaF} holds and $(\norm{w}/\la,\la)$ are sufficiently small. In fact, it is shown in~\cite{Liang14} that FB identifies in finite time the manifold of any non-degenerate minimizer $\xxs$. As a corollary, if condition~\eqref{eq-inj-etaF} holds at $x_0$ and $(\norm{w}/\la,\la)$ are sufficiently small, then we recover Theorem~\ref{thm-fb}. These results shed light on the typical convergence behavior of FB observed in such circumstances (e.g. in compressed sensing problems). 

\begin{remark}[Local linear convergence]
The FB generally exhibits a global sublinear $O(1/n)$ convergence rate in terms of the objective function. However, under partial smoothness of $J$, it is shown in~\cite{Liang14} that once the active manifold is identified, the FB algorithm enters a local linear convergence regime ($Q$-linear in general and $R$-linear if $\Mm$ is a linear manifold), whose rate can be characterized precisely in terms of the condition number of $\Phi_{T_{x_0}}$.
\end{remark}

\subsection{Related Works}
Finite support identification and local $R$-linear convergence of FB to solve~\eqref{eq-lagrangian} is established in \cite{bredies2008linear} under either a very restrictive injectivity assumption, or a non-degeneracy assumption that is a specialization of ours to the $\ell_1$ norm. A similar result is proved in \cite{hale07}. The $\ell_1$ norm is a partly smooth function and is therefore covered by Theorem~\ref{thm-fb}. \cite{negahban2010unified} proved $Q$-linear convergence of FB to solve~\eqref{eq-lagrangian} with a data fidelity satisfying restricted smoothness and strong convexity assumptions, and $J$ a so-called convex decomposable regularizer. Again, the latter falls within the class of partly smooth functions, and their result is then subsumed by our analysis. 

For general programs, a variety of algorithms, such as proximal and projected-gradient schemes were observed to have the finite identification property of the active manifold. In \cite{HareLewis04,hare2007identif}, the authors have shown finite identification of manifolds associated to partly smooth functions via the (sub)gradient projection method, Newton-like methods, and the proximal point algorithm. Their work extends that of e.g. \cite{Wright-IdentSurf} on identifiable surfaces from the smooth constrained convex case to a general non-smooth setting. Using these results, \cite{HareFB11} considered the algorithm \cite{Tseng09} to solve~\eqref{eq:minpbgen} when $h$ is $C^2$, $K=1$, $A_1=\Id$, and $g_1$ is simple and partly smooth, but not necessarily convex, and proved finite identification of the active manifold. However, the convergence rates remain an open problem in all these works.


\section{Summary and Perspectives}

In this chapter, we have reviewed work covering a large body of literature on the regularization of linear inverse problems. We also showed how these previous works can be all seen as particular instances of a unified framework, namely sensitivity analysis for minimization of convex partly smooth functions. We believe this general framework is the one that should be adopted as long as one is interested in studying fine properties and guarantees of these regularizers, and in particular when the stability of the low-complexity manifold associated to the data to recover is at stake.

This analysis is however only the tip of the iceberg, and there is actually a flurry of open problems to go beyond the theoretical results presented in this chapter. We list here a few ones that we believe are important avenues for future works:
\begin{itemize}
	\item \textit{Non-convexity and/or non-finiteness:} in this chapter, for the sake of simplicity, we focused on smooth convex fidelity terms and finite-valued convex regularizers. All the results stated in this chapter extend readily to proper lower semicontinuous convex regularizers, since any such a function is subdifferentially regular.
	Generalizations of some of the results to non-convex regularizers is possible as well, though some regularity assumptions are needed. This is of practical importance to deal with settings where $\Phi$ is not a linear operator, or to impose more agressive regularization (for instance when using $\ell^p$ functional with $0 \leq p  < 1$ instead of the $\ell^1$ norm). There are however many difficulties to tackle in this case. For instance, regularity properties that hold automatically for the convex case have to be either imposed or proved. Another major bottleneck is that some of the results presented here, if extended verbatim, will only assess the recovery of a stationary/critical point. The latter is not a local minimum in general, and even less global. 
	%
	\item \textit{Dictionary learning:} a related non-convex sensitivity analysis problem is to understand the recovery of the dictionary $D$ in synthesis regularization (as defined in Section~\ref{sec-synth-sparsity}) when solving problems of the form
	\eq{
		\umin{\{\al_k\}_k, D \in \Dd} \sum_k \frac{1}{2}\norm{y-\Phi D \al_k}^2 + \la J_0(\al_k)
	}
	where the $(y_k)_{k}$ are a set of input exemplars and $\Dd$ stands for the set of constraints imposed on the dictionary to avoid trivial solutions.
	Such a non-convex variational problem is popular to compute adapted dictionaries, in particular when $J_0=\norm{\cdot}_1$, see~\cite{elad2010sparse} and  references therein. 
	Although the dictionary learning problem has been extensively studied when $J_0=\norm{\cdot}_1$, most of the methods lack of theoretical guarantees. The theory of dictionary learning is only beginning to develop, see e.g.~\cite{GribonvalSchnass10,Jenatton,Spielman,Agarwal}.
	Tackling other regularizers, including analysis $\lun$ of the form $J=J_0 \circ D^*$ is even more difficult, see e.g.~\cite{PeyreFadiliAnalysis11,ChenPock12} for some computational schemes.
	\item \textit{Infinite dimensional problems:} we dealt in this chapter with finite-dimensional vector spaces. It is not straightforward to extend these results to infinite-di\-men\-sio\-nal cases. As far as $\ell_2$-stability is concerned, the constants involved in the upper bounds depend on the dimension $\N$, and the scaling might diverge as $N \rightarrow +\infty$. We refer to Section~\ref{sec-litterature-l2} for previous works on convergence rates of Tikhonov regularization in infinite-dimensional Hilbert or Banach spaces. Extending Theorem~\ref{thm-stability} for possibly non-reflexive Banach spaces is however still out of reach (non-reflexivity is a typical degeneracy when considering low-complexity regularization). 
	There exists however some extensions of classical stability results over spaces of measures, such as 
	weak convergence~\cite{bredies-inverse2013}, exact recovery~\cite{candes-towards2013,deCastro-exact2012} and stable support recovery~\cite{duval2013spike}.
	%
	\item \textit{Compressed sensing:} as highlighted in Sections~\ref{sec-cs-dual-certif} and~\ref{sec-cs-modelstab}, the general machinery of partly smooth regularizers (and the associated dual certificates) is well adapted to derive optimal recovery bounds for compressed sensing. Unfortunately, this analysis has been for now only applied to norms ($\norm{\cdot}_1$, $\norm{\cdot}_{1,\Bb}$, $\norm{\cdot}_*$ and $\norm{\cdot}_\infty$). Extending this framework for synthesis and analysis regularizers (see Sections~\ref{sec-synth-sparsity} and~\ref{sec-analysis-regularizers}) is a difficult open problem.   
	\item \textit{Convergence and acceleration of the optimization schemes:} Section~\ref{sec-fb} showed how partial smoothness can be used to achieve exact manifold identification after a finite number of iterations using the FB algorithm. This in turn implies a local linear convergence of the iterates, and raises the hope of acceleration using either first-order or second-order information fo the function along the identified manifold (in which we recall it is $C^2$). Studying such accelerations and their guarantees as well as extending this idea other proximal splitting schemes is thus of practical importance to tackle more complicated problems such as e.g. \eqref{eq:reg-noiseless}, \eqref{eq-constraint-epsilon} or \eqref{eq-constraint-gamma}. 
\end{itemize}

\section*{Acknowledgements} 

This work has been supported by the European Research Council (ERC project SIGMA-Vision).
We would like to thank our collaborators 
Charles Deledalle, 
Charles Dossal, 
Mohammad Golbabaee and
Vincent Duval
who have helped to build this unified view of the field.

\bibliographystyle{plain}
\bibliography{all}

\end{document}